\documentclass{svjour3}
\usepackage{amsmath,amsfonts,amssymb}

\def\xii{{\boldsymbol\xi}}
\def\etaa{{\boldsymbol\eta}}
\def\zetaa{{\boldsymbol\zeta}}
\def\phii{{\boldsymbol\varphi}}

\def\nuu{{\boldsymbol\nu}}
\def\uu{{\boldsymbol u}}
\def\thetaa{{\boldsymbol\theta}}

\renewcommand{\Pr}{{\mathbb P}}
\newcommand{\E}{{\mathbb E}}
\newcommand{\R}{{\mathbb R}}

\newcommand{\Z}{{\mathbb Z}}

\newcommand{\Sp}{{\mathbb S}}
\def\RO{\R^d\setminus\{\boldsymbol 0\}}
\def\Vfi{V_{\boldsymbol{\varphi}}}

\usepackage{color}
\definecolor{c20}{rgb}{0.,0.7,0.}
\definecolor{c30}{rgb}{0.,0.,1.}
\definecolor{c40}{rgb}{1,0.1,0.7}
\definecolor{c50}{rgb}{1,0,0}
\definecolor{c60}{rgb}{0,0.9,0.1}

\begin{document}

\title{Asymptotic Expansion of Gaussian Chaos via Probabilistic Approach}

\author{Enkelejd  Hashorva \and Dmitry Korshunov \and Vladimir I. Piterbarg}

\institute{E. Hashorva \at University of Lausanne, (HEC Lausanne), UNIL-Dorigny, 1015 Lausanne, Switzerland\\
              \email{Enkelejd.Hashorva@unil.ch}\\
           \and
           D. Korshunov \at Lancaster University, UK, and
           Sobolev Institute of Mathematics, Novosibirsk, Russia\\
              \email{d.korshunov@lancaster.ac.uk}\\
           \and
           V. I. Piterbarg, Lomonosov Moscow State University, Moscow, Russia\\
              \email{piter@mech.math.msu.su} 
}

\date{Received: date / Accepted: date}

\maketitle

\begin{abstract}
For a centered $d$-dimensional Gaussian random vector
$\xii =(\xi_1,\ldots,\xi_d)$ and a homogeneous function
$h:\R^d\to\R$ we derive asymptotic expansions
for the tail of the Gaussian chaos 
$h(\xii)$ given the function $h$ is sufficiently smooth. 
Three challenging instances of the Gaussian chaos are  
the determinant of a Gaussian matrix, the Gaussian orthogonal ensemble 
and the diameter of random Gaussian clouds.
Using a direct probabilistic asymptotic method, we investigate 
both the asymptotic behaviour of the tail distribution 
of $h(\xii)$ and its density at infinity and then discuss possible extensions 
for some general $\xii$ with polar representation.


\end{abstract}

\section{Introduction and Main Results}

Let $\xii=(\xi_1,\ldots,\xi_d)$ be a centered Gaussian random vector in $\R^d$,
$d\ge 2$, with covariance matrix $B$, $B_{ij}:=\E\xi_i\xi_j$.
Let $h:\R^d\to\R$ be a homogeneous function of order $\alpha>0$, that is,
$h(x\boldsymbol t)=x^\alpha h(\boldsymbol t)
$ for all $x>0$ and $\boldsymbol t=(t_1,\ldots,t_d) \in\R^d$.
Simple examples for $h$ are
$h(\boldsymbol t)= \prod_{i=1}^d |t_i|^{\gamma_i}$ of order
$\alpha=\gamma_1+\ldots+\gamma_d$ and
$h(\boldsymbol t)= \sum_{i=1}^d |t_i|^\alpha$ also of order $\alpha$.

We say that the random variable $h(\xii)$ is a {\it Gaussian chaos
of order} $\alpha$. In the literature, the term Gaussian chaos
of integer order $\alpha$ is traditionally reserved for the case
where $g$ is a homogeneous polynomial of degree $\alpha$---this
case goes back to Wiener (1938) where polynomial chaos
processes were first time introduced---and by this reason it is
spoken about as Wiener chaos in the Gaussian case.
Here we follow the extended version of the term Gaussian chaos.

This contribution is concerned with the asymptotic behavior
of the tail distribution of Gaussian chaos
$h(\xii)$ and its density at infinity.
We suppose that $h$ is not negative, that is, for some $\boldsymbol{x}$,
$h(\boldsymbol{x})>0$, otherwise our problem is trivial.
The important contributions in this area are
Hanson and Wright (1971) where an upper rough bound is obtained
for the tail of $h(\xii)$ in the case of polynomial $h$ of degree $2$,
Borell (1978, Theorem 2.2),
Arcones and Gin\'e (1993, Corollary 4.4),
Janson (1997, Theorem 6.12), and Lata{\l}a (1999, 2006)
where some lower and upper bounds are derived
in the case of polynomial $h$ of general degree $\alpha\ge 2$
(see also Lehec (2011)). A closely related study is
devoted to the derivation of lower and upper bounds for the distribution of
the multiple Wiener--It\^o integrals with respect to a white noise,
see Major (2005, 2007).
In all these papers the estimation
of the distribution tail of $h(\xii)$ is based on upper
bounds for the moments of $h(\xii)$; clearly this technique
cannot help with exact asymptotics for the tails.

In this contribution we shall focus on the Gaussian framework,
so the random vector $\xii$ introduced above is equal in
distribution to $\sqrt{B}\etaa$ where the coordinates of
$\etaa=(\eta_1,\ldots,\eta_d)$ are independent with
standard normal distribution.
Thus for any $x$ positive
\begin{equation}\label{eqaA}
\Pr\{h(\xii)>x\}=\Pr\{h(\sqrt{B}\etaa)>x\}=\Pr\{g(\etaa)>x\},
\end{equation}  
with  $g(\boldsymbol{u})=h(\sqrt{B}\boldsymbol{u})$.
The function $g:\R^d\to\R$ is homogeneous of order
$\alpha$ as $h$ is. In this standard way the problem for a general
covariance matrix may be reduced to that with identity matrix.

The distribution of $g(\etaa)$ may contain an atom at zero point, 
$\Pr\{g(\etaa)=0\}\ge 0$.
Remarkably, the distribution of $g(\etaa)$ restricted to $\R\setminus\{0\}$ 
always possesses a density function $p_{g(\etaa)}(x)$, $x\neq 0$, 
and the following representation
\begin{equation}\label{eqden}
p_{g(\etaa)}(x) =
\frac{1}{\alpha x}\bigl(\E\{\|\etaa\|^2;g(\etaa)>x\}
-d\cdot\Pr\{g(\etaa)>x\}\bigr)
\end{equation}
is valid for any $x>0$, see Lemma \ref{dens.via.tail} below.

Motivated by \eqref{eqaA} and \eqref{eqden},
in the following we shall formulate our results 
for the Gaussian chaos $g(\etaa)$.

For the Gaussian chaos, at least two approaches are available for the asymptotic
analysis of the tail distribution. The first approach is based on the asymptotic 
Laplace method and the second one exploits the rotation invariance of the standard 
normal distribution in $\R^d$ and may be regarded as a probabilistic approach. 
In the present paper we follow the probabilistic approach which is particularly
convenient for study of the elliptic chaos, see Theorem \ref{thPolar} below.
Earlier in our short note \cite{HKP}, we suggested to follow the asymptotic 
Laplace method in order to derive tail asymptotics for $g(\etaa)$. 
Notice that the Laplace method gives less information on the most probable
event where chaos large deviations occur. The advantage of the Laplace method
is that it is easily applicable to so-called Weibullian chaos;
the corresponding results will be presented in a forthcoming paper. 

So, in this contribution our analysis is based on the rotation invariance of 
the standard normal distribution. That is, for a $d$-dimensional centered Gaussian
random vector $\etaa$ with identity covariance matrix, the polar representation
\begin{equation}\label{eta.polar}
\etaa\overset{d}=\chi\zetaa
\end{equation}
holds in distribution where $\chi$ and $\zetaa$ are independent,
$\chi^2=\sum_{i=1}^d\eta_i^2$ has $\chi^2$-distribution with $d$
degrees of freedom and $\zetaa$ is uniformly distributed on the unit
sphere $\Sp_{d-1}\subset\R^d$. Hence by the homogeneity property of $h$
for any $x>0$ we have
\begin{equation}\label{g.chi.mu}
\Pr\{g(\etaa)>x\}=\Pr\{\chi^\alpha g(\zetaa)>x\}.
\end{equation}
In the sequel $g$ is assumed to be continuous, so that $g(\zetaa)$ 
is a non-negative bounded random variable.
The random variable $\chi^\alpha$ has the density function
\begin{equation}\label{chi.alpha}
p_{\chi^\alpha}(x) =
\frac{1}{\alpha2^{d/2-1}\Gamma(d/2)}x^{d/\alpha-1}e^{-x^{2/\alpha}/2},
\quad x>0,
\end{equation}
of Weibullian type with index $2/\alpha$ which is subexponential
density if $\alpha>2$, see e.g., Foss et al.\ \cite[Sect. 4.3]{FKZ}.
By this reason, the tail behaviour of the product
$\chi^\alpha g(\zetaa)$ heavily depends on the maximum
of the function $g$ on the unit sphere $\Sp_{d-1}$.
Denote
$$
\hat g:=\max_{\boldsymbol v\in\Sp_{d-1}}g(\boldsymbol v)
\quad\mbox{and}\quad \mathcal M:=\{\boldsymbol v\in\Sp_{d-1}:g(\boldsymbol v)=\hat g\}.
$$
We shall consider two different cases of the structure of the
set $\mathcal M$:

\begin{description}
\item[(i)] $\mathcal M$ consists of a finite number of isolated points.

\item[(ii)] $\mathcal M$ is a sufficiently smooth manifold of
positive dimension $m$, $1\le m\le d-2$, on the unit sphere.
\end{description}

In the second case we assume that $\mathcal M$ has no boundary
which particularly assumes that $m\neq d-1$. This restriction
comes from the observation that the existence of a boundary
of the set of the points of maximum $\mathcal M$ strongly
contradicts the condition that the function $g$ is at least
twice continuously differentiable with non-degenerate approaching
of its maximum.

\subsection{The case of finite $\mathcal M$}

Here we consider a homogeneous continuous function $g:\R^d\to\R$
of order $\alpha>0$ such that $\mathcal M$ consists of
a finite number of points, say
$$
\mathcal M=\{\boldsymbol v_1,\ldots,\boldsymbol v_k\}.
$$

Let in the following $g\in C^2(\RO)$.
For every point $\boldsymbol v\in\Sp_{d-1}$, 
denote by $g''_{d-1}(\boldsymbol v)$ a Hessian matrix at point $\boldsymbol v$ 
of the function $g$ restricted to the hyperplane 
tangent to the sphere $\Sp_{d-1}$ at point $\boldsymbol v$, 
that is, restricted to the hyperplane $\boldsymbol v+\mathcal L$ 
where $\mathcal L=\{\boldsymbol u\in\R^d:(\boldsymbol u,\boldsymbol v)=0\}$.
More precisely, we fix an orthogonal system of vectors in $\mathcal L$,
say $\uu_1$, \ldots, $\uu_{d-1}$ and consider the function 
$g_{d-1}(t_1,\ldots,t_{d-1}):=g(\boldsymbol v+t_1\uu_1+\ldots+t_{d-1}\uu_{d-1})$
whose Hessian matrix is denoted by $g''_{d-1}(\boldsymbol v)$.

Assume that  for every $j=1$, \ldots, $k$
\begin{eqnarray}\label{non.gener.f}
\det\Bigl(\frac{g''_{d-1}(\boldsymbol v_j)}{\alpha\hat g}
-I_{d-1}\Bigr) &<& 0,
\end{eqnarray}
where $I_n$ stands for the identity matrix of size $n$.
As follows from Lemma \ref{l:car.local},
$$
g''_{d-1}(\boldsymbol v_j)-\hat g\alpha I_{d-1}
$$
is just a Hessian matrix of the function $g$ along the unit sphere at point 
$\boldsymbol v_j$; the latter is explained 
in more detail after Theorem \ref{th:m=0}. 
The condition \eqref{non.gener.f} says that
$\boldsymbol v_j\in\Sp_{d-1}$ is the point of non-degenerate maximum
of the function $g$ restricted to $\Sp_{d-1}$.
Then the following result holds.

\begin{theorem}\label{th:m=0}
Let $g\in C^{2r+2}( \RO )$ for some $r\ge 0$.
Then the following asymptotical expansion takes place, as $x\to\infty$:
\begin{equation}\label{AE1}
\Pr\{g(\etaa)>x\}=(x/\hat g)^{-1/\alpha}
e^{-(x/\hat g)^{2/\alpha}/2}
\biggl(h_0+\sum_{i=1}^r h_ix^{-2i/\alpha}+o(x^{-2r/\alpha})\biggr),
\end{equation}
where coefficients $h_0$, \ldots, $h_r\in\R$ only depend on $\alpha$,
$\hat g$, and derivatives of $g(\phii)$ at points $\phii_j$
{\rm(}for the definition of $g(\phii)$ see below after \eqref{spheric}{\rm)}; 
in particular,
\begin{eqnarray}\label{0h0.1}
h_0 &:=& \frac{1}{\sqrt{2\pi}}
\sum_{j=1}^k
\Bigl|\det\Bigl(\frac{g''_{d-1}(\boldsymbol v_j)}{\alpha\hat g}
-I_{d-1}\Bigr)\Bigr|^{-1/2}
\end{eqnarray}
Moreover, the density function of $g(\etaa)$
satisfies the following relation, as $x\to\infty$:
\begin{eqnarray}\label{cor:m=0}
p_{g(\etaa)}(x)=(x/\hat g)^{1/\alpha-1}
e^{-(x/\hat g)^{2/\alpha}/2} \biggl(\frac{h_0}{\alpha\hat g}
+\sum_{i=1}^r \tilde h_ix^{-2i/\alpha}+o(x^{-2r/\alpha})\biggr).
\end{eqnarray}
\end{theorem}

Notice that if $\alpha=2$, then the tail distribution
of $g(\etaa)$ is asymptotically proportional to its
density function $p_{g(\etaa)}(x)$ as $x\to\infty$
with multiplier $2\hat g$.

It is essential assumption that the function $g$ is at least in $C^2$ and 
that its Hessian along the unit sphere is non-degenerate on $\mathcal M$.
If it is not so, then the tail asymptotics may be quite specific and requires 
additional investigation---especially in the case where $g\not\in C^2$. 
An example of natural Gaussian chaos with degenerate Hessian 
along the unit sphere is discussed below, 
see Example 9 in Section \ref{sec:discussion}.

Sometimes it is more convenient to pass to some local
coordinates on the sphere. Let $V_j\subset\R^d$
be a neighborhood of the point $\boldsymbol v_j\in\mathcal M$ and let $h_j$ be a twice differentiable bijection from the open cube
$(0,2)^d$ to $V_j$ such that $h_j$ is a bijection from
$\{\boldsymbol z\in(0,2)^d:z_d=1\}$ to $V_j\cap\Sp_{d-1}$.
Denote by
\[
(g\circ h_j)_{d-1}''(\boldsymbol z):=\left[\frac{\partial^2(g\circ h_j)(\boldsymbol z)}
{\partial z_i\partial z_l}\right]_{i,l=1,\ldots,d-1},
\quad \boldsymbol z\in (0,1)^{d-1},
\]
the Hessian matrix of $g\circ h_j$ restricted to the first $d-1$
coordinates---it is a $(d-1)\times(d-1)$ matrix---and 
write  $\boldsymbol z_j\in (0,2)^{d-1}\times\{1\}$
for a point satisfying  $h_j(\boldsymbol z_j)=\boldsymbol v_j$.
We will prove in Lemma \ref{l:car.local} that, at every point
$\boldsymbol v_j\in\mathcal M$, the following equality holds:
\begin{eqnarray}\label{cart.local}
\det(g''_{d-1}(\boldsymbol v_j)-(\alpha\hat g)I_{d-1})
&=& \frac{\det (g\circ h_j)_{d-1}''(\boldsymbol z_j)}
{(\det J_j(\boldsymbol z_j))^2},
\end{eqnarray}
where $J_j(\boldsymbol u)$ is the Jacobian matrix of $h_j$.
Then the representation \eqref{0h0.1} for the constant $h_0$
can be rewritten in terms of local coordinates as follows:
\begin{eqnarray}\label{0h0.2}
h_0 &:=& \frac{1}{\sqrt{2\pi}} (\alpha\hat g)^{\frac{d-1}{2}}
\sum_{j=1}^k \frac{|\det J_j(\boldsymbol z_j)|}
{\sqrt{|\det (g\circ h_j)_{d-1}''(\boldsymbol z_j)|}}.
\end{eqnarray}

A particular example is given by the hyperspherical
coordinates, $\boldsymbol v=(r,\phii)$, with Jacobian
\[
\det J(r,\phii)=r^{d-1}\sin^{d-2}\varphi_1\ldots\sin\varphi_{d-2}
=r^{d-1}\det J(1,\phii),
\]
where $J(r,\phii)$ stands for the Jacobian matrix, that is,
\begin{align}\label{spheric}
& v_1=r\cos\varphi_1\nonumber\\
& v_2=r\sin\varphi_1\cos\varphi_2\nonumber\\
& \ldots\nonumber\\
& v_{d-1}=r\sin\varphi_1\sin\varphi_2\ldots\sin\varphi_{d-2}
\cos\varphi_{d-1}\nonumber\\
& v_d=r\sin\varphi_1\sin\varphi_2\ldots\sin\varphi_{d-2}\sin\varphi_{d-1},
\end{align}
with  $\phii=(\varphi_1,\ldots,\varphi_{d-1})
\in\Pi_{d-1}:=[0,\pi)^{d-2}\times\lbrack0,2\pi)$ the angular coordinates 
of $\boldsymbol v$, $r=\|\boldsymbol v\|$.
As usual, the topology in the set $\Pi_{d-1}$ is induced
by the topology on the unit sphere, in particular,
all points of $\Pi_{d-1}$ are inner points.
Changing in such a way variables, we have
(we set $g(\phii)=g(\boldsymbol v/\|\boldsymbol v\|)$;
the function $g(\phii)$ is continuous too;
hereinafter we denote by the same symbol $g$ two formally different functions,
on $\R^d$ and on $\Pi_{d-1}$, but this hopefully does not lead to any confusion)
$$
\hat g=\max_{\phii\in\Pi_{d-1}}g(\phii)
\quad\mbox{and}\quad
{\mathcal M}_\varphi:=\{\phii\in\Pi_{d-1}:g(\phii)=\hat g\}.
$$
Denote by $g''(\phii)$ the Hessian matrix of
$g(\varphi_1,\ldots,\varphi_{d-1})$. So, in the particular case of
hyperspherical coordinates, the equality \eqref{0h0.2} implies that
\begin{eqnarray}\label{0h0.3}
h_0 &:=& \frac{1}{\sqrt{2\pi}} (\alpha\hat g)^{\frac{d-1}{2}}
\sum_{j=1}^k \frac{|\det J(1,\phii_j)|}{\sqrt{|\det g''(\phii_j)|}}.
\end{eqnarray}
Then the condition \eqref{non.gener.f} requires that
the maximum of the function $g(\phii)$, $\phii\in\Pi_{d-1}$, 
is non-degenerate at points $\phii_1, \ldots, \phii_k$.

Some results for the case of isolated points of maximum may be
also found in Breitung and Richter \cite{BR1996}. 
It seems that the proofs of the asymptotic expansions
of Theorems 4 and 5 in \cite{BR1996} cannot be considered as self-contained. 
For instance, it is only proven in \cite[Lemma 3]{BR1996} that $\gamma_2=0$. 
In order to prove their Theorems 4 and 5 as they are stated, 
it is necessary to prove that all $\gamma_{2m=0}$.
In addition, at the beginning of their proofs of Theorems 4 and 5 
Breitung and Richter \cite{BR1996} write that all coefficients $a_i$ 
for odd $i$ are zero because of Theorem 4.5 
in Fedoryuk \cite[p. 82]{Fedoryuk}. This argument does not help
because Theorem 4.5 in Fedoryuk \cite[p. 82]{Fedoryuk} 
is not about asymptotic expansion of $F(A;(1+z)^{1/2})$ what is required
by the authors, so that it is irrelevant to the issue considered.

\subsection{The case of a manifold}

Now consider the case where $\mathcal M\subset\Sp_{d-1}$ is for some 
$m\in\{1,\ldots,d-2\}$ an $m$-dimensional manifold of finite volume 
and has no boundary.

Fix some $r\in\Z^+$. We assume that the manifold $\mathcal M$
is $C^{2r+2}$-smooth. 

We suppose that the rank of the matrix
$A_{d-1}(\boldsymbol v):=
\frac{g''_{d-1}(\boldsymbol v)}{\hat g\alpha}-I_{d-1}$ of
size $d-1$ is equal to $d-1-m$ for every
$\boldsymbol v\in\mathcal M$. Denote by
$\det\bigl(\frac{g_{d-1-m}''(\boldsymbol v)}{\hat g\alpha}-I_{d-1-m}\bigr)$
any non-zero $(d-1-m)$-minor of the matrix $A_{d-1}(\boldsymbol v)$;
notice that all $(d-1-m)$-minors are equal one to another,
by using orthogonal transform and set
\begin{eqnarray}\label{dh0.1}
h_0 &:=& \frac{1}{(2\pi)^{\frac{m+1}{2}}} \int_{\mathcal M}
\Bigl|\det\Bigl(\frac{g_{d-1-m}''(\boldsymbol v)}
{\hat g\alpha}-I_{d-1-m}\Bigr)\Bigr|^{-1/2}{\rm d}V,
\end{eqnarray}
where $dV$ is the volume element of
$\mathcal M\subset\Sp_{d-1}$.

\begin{theorem}\label{th:m=1.d-2}
Assume that the above conditions on $\mathcal M$ are fulfilled
and that $g\in C^{2r+2}( \RO )$.
Then the following asymptotical expansion takes place, as $x\to\infty$:
\begin{equation}\label{AE2}
\Pr\{g(\etaa)>x\}=(x/\hat g)^{\frac{m-1}{\alpha}}
e^{-(x/\hat g)^{2/\alpha}/2}
\biggl(h_0+\sum_{i=1}^r h_ix^{-2i/\alpha}+o(x^{-2r/\alpha})\biggr),
\end{equation}
where coefficients $h_1$, \ldots, $h_r\in\R$ only depend on $\alpha$,
$\hat g$, and derivatives of $g(\phii)$ on $\mathcal M_\varphi$.

Moreover, the  density function of $g(\etaa)$
satisfies the following relation, as $x\to\infty$:
\begin{eqnarray}\label{cor:m=1.d-2}
p_{g(\etaa)}(x) &=& (x/\hat g)^{\frac{m+1}{\alpha}-1}e^{-(x/\hat g)^{2/\alpha}/2}
\biggl(\frac{h_0}{\alpha\hat g}
+\sum_{i=1}^r\tilde h_ix^{-2i/\alpha}+o(x^{-2r/\alpha})\biggr).
\end{eqnarray}
\end{theorem}

Notice that if the manifold $\mathcal M$ has boundary points,
then asymptotic expansion becomes more complicated. 
In general, boundary points have no impact on the leading constant $h_0$.
Boundary makes strong contribution on further terms. For instance,
if $d=3$ and $\mathcal M$ is a line-segment on the unit sphere $\Sp_2$,
then the term $x^{-1/\alpha}$ appears in the parentheses 
of the expansions \eqref{AE2} and \eqref{cor:m=1.d-2}; 
the corresponding calculations in the neighborhood of the boundary 
may be rather specific compared to those in Lemmas \ref{l:1.f-1} and 
\ref{l:1.f-d} below. The main reason for this comes from the fact that,
in most cases, the function $g$ is not in $C^2$ on the boundary;
this is clearly demonstrated by the function 
$g(x)=-x^2$ for $x\le 0$ and $g(x)=0$ for $x\ge 0$.

By the same reasons as in the case of finite $\mathcal M$
we have the following representation for the constant $h_0$
in terms of the spherical coordinates:
\begin{eqnarray}\label{dh0.2}
h_0 &:=& \frac{1}{(2\pi)^{\frac{m+1}{2}}}
(\alpha\hat g)^{\frac{d-1-m}{2}}
\int_{\mathcal M_\varphi} \frac{|\det J(1,\phii)|}
{\sqrt{|\det g_{d-1-m}''(\phii)|}}{\rm d}\Vfi,
\end{eqnarray}
where ${\rm d}\Vfi$ is the volume element of
$\mathcal M_\varphi\subset\Pi_{d-1}$.

The organisation of the rest of the paper is as follows.
In Section \ref{sec:discussion} we discuss our main results
and provide several examples that concern different
cases for the dimension of $\mathcal M$.
Proofs of the main results are presented in Sections
\ref{sec:proof.0} and \ref{sec:proof.d}.

In \cite{HKP}, a preliminary version of Theorem \ref{th:m=1.d-2} was announced.
Precisely, the relation 
$$
\Pr\{g(\etaa)>x\}=(x/\hat g)^{\frac{m-1}{\alpha}}
e^{-(x/\hat g)^{2/\alpha}/2}\bigl(h_0+O(x^{-2/\alpha})\bigr)
\quad\mbox{ as }x\to\infty
$$ 
was stated without proof under the assumption that the function $g$ 
is three times differentiable. It was suggested to follow the asymptotic 
Laplace method in order to derive this relation.
One of the goals of the present paper is to provide a self-contained 
geometric proof of asymptotic expansion with $r+1$ 
terms under correct smoothness conditions.

\subsection{Elliptical chaos}

Before presenting several examples we show how our results can be extended 
for elliptical chaos or more generally for the chaos of polar random vectors. 
Consider therefore in the following $\xii$ such that \eqref{rF} holds with 
$\chi>0$ some random variable being independent of $\zetaa$. 
Crucial properties used in the Gaussian case are {\bf a)} $\chi$ has 
distribution function in the Gumbel max-domain of attraction, and {\bf b)} 
the random vector $\zetaa$ has a $d-1$ dimensional subvector which possesses 
a positive density function. The first property a) means that for any 
$t\in\R$
\begin{equation}\label{rG}
\Pr\{\chi>x+t/w(x)\} \sim e^{-t} \Pr\{\chi>x\}
\quad\mbox{as }x\uparrow x_+,
\end{equation}
with $w$ a positive scaling function and $x_+$ 
the upper endpoint of the distribution function of $\chi$ 
(in the Gaussian case $w(x)= x,$ and $x_+=\infty$). 
We abbreviate \eqref{rG} as  $\chi \in GMDA(w,x_+)$. 
Condition \eqref{rG} is satisfied by a large class of 
random variables, for instance if $\chi$ is such that
$$
\Pr\{\chi>x\} \sim c_1 x^a e^{-c_2x^\beta}
\quad\mbox{as } x\to\infty
$$
for some $c_1$, $c_2$, $\beta$ positive and $a\in\R$, then \eqref{rG}
holds with $w(x)= \beta c_2 x^{\beta-1}$. Notice that $c \chi^\beta $ 
is in the Gumbel max-domain of attraction for any $c>0$ and $\beta>0$ 
if and only if $\chi$ is in the Gumbel max-domain of attraction.

In order to relax the assumption on $\zetaa$ note first that
in hyperspherical coordinates the random vector
$\zetaa=(\zeta_1,\ldots,\zeta_d)\in\Sp_{d-1}$ can be written as
$\nuu=(\nu_1,\ldots,\nu_{d-1})\in\Pi_{d-1}$.
Since $\zetaa$ is uniformly distributed on the unit sphere $\Sp_{d-1}$
in $\R^d$, the density function of the random vector
$\nuu=(\nu_1,\ldots,\nu_{d-1})\in\Pi_{d-1}$
equals $\frac{|\det J(1,\phii)|}{{\rm mes}\,\Sp_{d-1}}$, $\phii\in\Pi_{d-1}$.

If $\chi$ is some positive random variable, then $\xii$ is an 
elliptically symmetric random vector. 
When $\chi^2$ is chi-square distributed with $d$ degrees of freedom
we recover as special case of elliptically symmetric random vectors
the Gaussian ones. In particular for the Gaussian case we have
\begin{equation}\label{chi}
\Pr\{\chi^\alpha> x\} \sim
\frac{1}{2^{d/2-1}\Gamma(d/2)} x^{(d-2)/\alpha}e^{-x^{2/\alpha}/2}
\quad\mbox{as } x\to\infty
\end{equation}
implying that $\chi^\alpha$ is in the Gumbel MDA with scaling function
$w(x)=x^{2/\alpha-1}/\alpha$, $x>0$. 
If we relax our assumption on
the distribution function of $\chi$ and simply assume \eqref{chi}
our previous results cannot be immediately re-formulated since the Gaussianity does not hold anymore.
It turns out that even the larger class of elliptically symmetric random vectors
for which $\chi$ satisfies \eqref{chi} is a strong (unnecessary)
restriction for the derivation of the tail asymptotics of $g(\etaa)$.
Indeed, we shall drop in the following the explicit distributional
assumption on $\zetaa$ assuming only that $\nuu$ possesses a 
positive bounded continuous density function, say $p_\nu(\phii)$.
Next, we present the counterpart of Theorem \ref{th:m=1.d-2},
i.e., as therein we shall impose the same conditions on $\mathcal{M}$.

\begin{theorem}\label{thPolar}
Assume that $g\in C^2( \RO )$.
If $\zetaa$ is such that $\chi^\alpha\in GMDA(w,x_+)$
and further the random vector $\nuu$ has a positive 
bounded continuous density function $p_\nu(\phii)$, then
\begin{eqnarray}
\Pr\{g(\etaa)>x\} &\sim&  
\frac{h_0}{(xw(x/\hat g))^{\frac{d-1-m}{2}}}\Pr\{\chi^\alpha>x/\hat g\}
\end{eqnarray}
as $x\uparrow \hat g x_+$, where
\begin{eqnarray}\label{ePo2}
h_0=(2\pi\hat g^2)^{\frac{d-1-m}{2}} \int_{\mathcal M_\varphi}
\frac{p_\nu(\phii)}{\sqrt{|\det g_{d-1-m}''(\phii)|}}{\rm d}\Vfi \in (0,\infty).
\end{eqnarray}
In particular $g(\etaa)\in GMDA(w,\hat g x_+)$.
\end{theorem}

\newcommand{\COM}[1]{}
\COM{{\bf Remark.} Let us consider two random chaoses with the same $g$
and $p_\nu(\phii)$ but with different $\chi_1$ and $\chi_2$.
Then it follows from Theorem \ref{thPolar} that, as $x\to\infty$, 
\begin{eqnarray*}
\Pr\{g(\etaa_2)>x\} &\sim&  
\Bigl(\frac{w_1(x/\hat g)}{w_2(x/\hat g)}\Bigr)^{\frac{d-1-m}{2}}
\frac{\Pr\{\chi_2^\alpha>x/\hat g\}}{\Pr\{\chi_1^\alpha>x/\hat g\}}
\Pr\{g(\etaa_1)>x\}.
\end{eqnarray*}
If $\etaa_1$ is standard normal vector, 
then $\chi_1^2$ is a chi-square random variable with $d$ degrees of freedom
satisfying \eqref{chi}, $w_1(x)=\frac{1}{\alpha}x^{\frac{2-\alpha}{\alpha}}$ 
and hence
\begin{eqnarray}\label{ratio.12}
\frac{\Pr\{g(\etaa_2)>x\}}{\Pr\{g(\etaa_1)>x\}} &\sim&  
\frac{\Pr\{\chi_2^\alpha>x/\hat g\}}{w_2^{\frac{d-1-m}{2}}(x/\hat g)}
(x/\hat g)^{\frac{2-\alpha}{\alpha}\frac{d-1-m}{2}-\frac{d-2}{\alpha}}
\frac{2^{d/2-1}\Gamma(d/2)}{\alpha^{\frac{d-1-m}{2}}} e^{(x/\hat g)^{2/\alpha}/2}.\nonumber\\[-2mm]
\end{eqnarray}

If $\nuu$ is elliptically symmetric, then $ h_0$ above 
can be directly calculated if the tail asymptotics of $g(\etaa)$ 
for the case that $\chi^2$ is a chi-square random variable 
with $d$ degrees of freedom is of the form 
$$ 
\Pr\{g(\etaa)>x\} \sim c x^{\frac{m-1}{\alpha}} e^{-(x/\hat g)^{2/\alpha}/2} 
\quad \text{ as } x\to\infty.
$$
Indeed, since $p_\nu(\phii)=J(1,\phii)/{\rm mes}\,\Sp_{d-1}$,  
our results above show that  $h_0$ in \eqref{ePo2} is simply given by 
\begin{eqnarray}\label{ePo2Pa}
h_0 = c  \hat g^{ \frac{m-1}{\alpha}}. 
\end{eqnarray}
Similarly, if as $x\to \infty$
$$ 
p_{g(\etaa)}(x) \sim c x^{(m+1)/\alpha -1} e^{-(x/\hat g)^{2/\alpha}/2},
$$
then we find that 
\begin{eqnarray}\label{ePo2Pb}
h_0 = \frac{c}{b L}  \hat g^{ \frac{m+1}{\alpha}}. 
\end{eqnarray}

Indeed by our results it follows that if $w(x)= x^{2/\alpha-1} /\alpha$
\begin{eqnarray*}
\Pr\{g(\etaa)>  x\}
&\sim& h_0 2^{d/2-1} \Gamma(d/2) 
\bigl(\alpha x/\hat g w(x/\hat g)\bigr)^{\frac{m+1-d}{2}}\Pr\{\chi^2>(x/\hat g)^{1/\alpha}\}\\
&\sim& h_0 2^{d/2-1} \Gamma(d/2) 
\bigl(x/\hat g (x/\hat g)^{2/\alpha -1} \bigr)^{\frac{m+1-d}{2}} \frac{2^{-d/2+1}}{\Gamma(d/2)}
(x/\hat g)^{2/\alpha (d/2-1)} e^{- (x/ \hat g)^{2/\alpha}/2}\\
&\sim& h_0 2^{d/2} (x/\hat g)^{2/\alpha \frac{m+1-d}{2}} (x/\hat g)^{2/\alpha (d/2-1)} e^{- (x/ \hat g)^{2/\alpha}/2}\\
&\sim& h_0 2^{d/2} (x/\hat g)^{ \frac{m+1-d}{\alpha} +d/\alpha - 2/\alpha} e^{- (x/ \hat g)^{2/\alpha}/2}\\
&\sim& h_0 (x/\hat g)^{ \frac{m-1}{\alpha}} e^{- (x/ \hat g)^{2/\alpha}/2}\\
\end{eqnarray*}
}

Higher order asymptotic expansions for $\Pr\{g(\etaa)>x\}$---when
$g\in C^{2r+2}$ and $\Pr\{\chi^\alpha>x\}$ possesses a suitable 
expansion---can also be derived; the asymptotics of the density
of $g(\etaa)$ can be given in terms of $\Pr\{g(\etaa)>x\}$ and
the scaling function $w$ if further $\chi^\alpha$ possesses a bounded density function
such that $\Pr\{\chi^\alpha>x\}\sim p_{\chi^\alpha}(x)/w(x)$ as $x\to x_+$. 
We shall omit these results here.
Additionally the tail asymptotics of $g(\etaa)$ can be found also
(using similar arguments as in the proofs of Theorem 3) if instead
of $\chi$ in the Gumbel max-domain of attraction we assume that
$\chi$ is in the Weibull max-domain of attraction, i.e.,
$$
\Pr\{\chi>x_+- s/x\} \sim s^\gamma \Pr\{\chi>x_+- 1/x\}
\quad\mbox{as }x\to\infty,
$$
for any $s>0$ with $\gamma\ge 0$ some given constant.
As in the Gumbel max-domain of attraction case, $\chi^\beta$ is
in the Weibull max-domain of attraction for some $\beta>0$
if and only if $\chi$ is in the Weibull max-domain of attraction,
see Resnick (1987) for details on max-domain of attractions.

\section{Discussion and Examples}
\label{sec:discussion}

In view of Theorems \ref{th:m=0} and \ref{th:m=1.d-2},
the Gaussian chaos is a subexponential random variable if
$\alpha>2$ (under the assumptions therein). Subexponentiality of
random variables is an important concept with various applications,
see e.g., Embrechts et al.\ \cite{EMB}, Foss et al.\ \cite{FKZ}. 
It is possible to show subexponentiality of Gaussian chaos under 
some weak conditions on the homogeneous function $h$.
As follows from the polar representation \eqref{eta.polar},
the Gaussian random vector
\begin{equation}\label{rF}
\xii=\sqrt{B}\etaa \stackrel{d}{=}\chi\sqrt{B}\zetaa
\end{equation}
 has covariance matrix $B$. Hence by the homogeneity property of $h$
for any $x>0$ we have
\begin{equation}\label{h.chi.mu}
\Pr\{h(\xii)>x\}
=\Pr\{\chi^\alpha h(\sqrt{B}\zetaa)>x\}.
\end{equation}
Assuming that $h(\sqrt{B}\zetaa)$ is a positive bounded
random variable, then in view of Cline and Samorodnitsky
\cite[Corollary 2.5]{Cline1994} the random variable
$h(\xii)$ is subexponential if $\alpha>2$ because
then $\chi^\alpha$ has the density function \eqref{chi.alpha}
which is of Weibullian type with index $2/\alpha<1$ and
subexponential by this reason.

As follows from the representation \eqref{h.chi.mu}, for $h$
bounded on $\Sp_{d-1}$, that is,
$\hat h:=\max\{h(\boldsymbol{u}):\|\boldsymbol{u}\|=1\}<\infty$,
\begin{align}\label{bound.uni}
\Pr\{h(\xii)>x\} & \le\Pr\{\chi^\alpha >x/\hat h\}\nonumber\\
& \le\frac{1}{\alpha2^{d/2-1}\Gamma(d/2)}
\int_{x/\hat h}^\infty y^{d/\alpha-1}e^{-y^{2/\alpha}/2}{\rm d}y.
\end{align}
This upper bound is explicit and valid for any homogeneous function $h$ 
as determined above. For the special case of {\it decoupled polynomial chaos}, 
upper bounds are known that are universal in $d$ but less explicit in $x$, 
see e.g., Lata{\l}a \cite[Corollary 1]{Latala2006}. 
See also upper bound by Arcones and Gin\'e \cite[Theorem 4.3]{Gine}
where the case of general polynomial chaos is considered;
the corresponding upper bound is not always better in $d$ 
than \eqref{bound.uni} and it is less explicit.

If (with necessity discontinuous) a function $h$ is unbounded on
the unit sphere $\Sp_{d-1}$, then it is possible that
$\Pr\{h(\sqrt{B}\zetaa)>x\}>0$ for any $x$.
Two cases of interest, which are also simple to deal with are
$h(\sqrt{B}\zetaa)$ is regularly varying with index
$\gamma>0$ and $h(\sqrt{B}\zetaa)$ has a Weibullian tail.

In the first case where the tail of $h(\sqrt{B}\zetaa)$ is heavier than that of $\chi$, by Breiman's theorem (see \cite{Breiman})
\[
\Pr\{h(\xii)>x\}\sim\E\{\chi^{\gamma\alpha}\}
\Pr\{h(\sqrt{B}\zetaa)>x\}\quad\text{as }x\to\infty.
\]
Also, in view of Jacobsen et al. (2009) the converse of the above holds, i.e.,
if $h(\xii)$ is a regularly varying random variable,
then $h(\sqrt{B}\zetaa)$ is regularly varying too,
with the same index. The second case that $h(\sqrt{B}\zetaa)$ 
has a Weibullian tail can be handled by applying
Lemma 3.2 in Arendarczyk and D\c{e}bicki (2011).

Notice that if the Gaussian vector $\xii$ with
covariance matrix $B$ has a singular distribution,
so that $\det B=0$, then $\xii$ is valued in the linear subspace 
$\mathcal{L}:=\{\sqrt{B}\boldsymbol{u}:\boldsymbol{u}\in\R^d\}$ 
of lower dimension $d^*<d$. Therefore, it is necessary to proceed 
to a Gaussian random vector $\xii^*$ in $\mathcal{L}$ of dimension $d^*$ 
and to a new function $h^*(u):=h(u)$ defined on $\mathcal{L}$. 
In this way the problem is reduced to that with non-degenerate Gaussian 
distribution. For example, let $d=2$ and $h(u_1,u_2)=u_2^{4}-u_1^{4}/2$. 
Let further $\eta_1$ and $\eta_2$ be two independent $N(0,1)$ 
random variables, and set $\xii_1=(\eta_1,\eta_2)$ and $\xii_2=(\eta_2,\eta_2)$. 
Then the tail of $h(\xii_1)=h(\eta_1,\eta_2)=\eta_2^{4}-\eta_1^{4}/2$
is equivalent to that of $\eta_2^{4}$ which is much heavier than the tail of
$h(\xii_2)=h(\eta_2,\eta_2)=\eta_2^{4}/2$.

\textbf{Example 1.}
Consider the chaos
$g(\etaa)=|\eta_1|^\alpha +\ldots+|\eta_d|^\alpha$ 
of order $\alpha>0$.

If $\alpha=2$, so that we deal with $\chi^2$-distribution
with $d$ degrees of freedom, then $\hat g=1$ and $\mathcal M$ is
the whole unit sphere $\Sp_{d-1}$ (which is a manifold of dimension
$d-1$) and as known, the density function of $g(\etaa)$ equals
$\frac{1}{2^{d/2}\Gamma(d/2)}x^{d/2-1}e^{-x/2}$.

If $\alpha<2$, then the function $y_1^{\alpha/2}+\ldots+y_d^{\alpha/2}$
is concave and, therefore, its maximum on the set
$y_1+\ldots+y_d=1$, $y_i\ge 0$,
is attained at the point $y_1=\ldots=y_d=1/d$.
Hence, $\hat g=d^{1-\alpha/2}$ and $\mathcal M$ consists of
$2^d$ points $(\pm1/\sqrt d,\ldots,\pm1/\sqrt d)$
(which is a manifold of zero dimension). Then, by \eqref{cor:m=0},
the density function of $g(\etaa)$ is equal to
$$
c_2x^{1/\alpha-1}e^{-d^{1-2/\alpha}x^{2/\alpha}/2}(1+O(x^{-2/\alpha}))
\quad\mbox{ as }x\to\infty,
$$
where 
\begin{eqnarray*}
c_2&=& 2^d\frac{1}{\alpha\hat g^{1/\alpha}}
\frac{1}{\sqrt{2\pi}\sqrt{\Bigl|\det\Bigl(\frac{g''_{d-1}(\boldsymbol v_j)}{\alpha\hat g}
-I_{d-1}\Bigr)\Bigr|}},
\end{eqnarray*}
where $g''_{d-1}(\boldsymbol v_j)=\frac{\alpha(\alpha-1)}{d^{\alpha/2-1}}I_{d-1}$, 
$j\le 2^d$. 
Therefore,
\begin{eqnarray*}
c_2&=& \frac{2^d}{\alpha d^{1/\alpha-1/2}
\sqrt{2\pi}(2-\alpha)^{(d-1)/2}};
\end{eqnarray*}
this is a special case of the result by Rootz\'en (1987, see (6.1));
see also Theorem 1.1 and Example 1.3 in Balkema et al. (1993).
We see that the power of $x$ in front of the exponent is changing
together with the dimension of the manifold $\mathcal M$.

If $\alpha>2$, then $\hat g=1$ and $\mathcal M$ consists
of $2d$ points $(0,\ldots,0,\pm1,0,\ldots,0)$ (which is again a
manifold of zero dimension) and, by Theorem \ref{th:m=0}
with $r=1$, the density function of $g(\etaa)$ is equal to
$$
c_3x^{1/\alpha-1}e^{-x^{2/\alpha}/2}(1+O(x^{-2/\alpha}))
\quad\mbox{ as }x\to\infty.
$$
Here the matrix $g''_{d-1}(\boldsymbol v_j)$ is zero at every
point $\boldsymbol v_j$ implying
$$
c_3=\frac{h_0}{\alpha}=\frac{2d}{\alpha\sqrt{2\pi}}.
$$
In this case $|\eta_i|^\alpha$ has subexponential density
$$
\frac{2}{\alpha\sqrt{2\pi}}x^{1/\alpha-1}e^{-x^{2/\alpha}/2},
\quad x> 0,
$$
and the density function of $g(\etaa)$ is asymptotically equivalent to $d$
multiple of the density function of $|\eta_i|^\alpha$ as
$x\to\infty$, see Foss et al. (2011, Chapter 4).
Here the observation that $\mathcal M$ consists of $d$ points
$(0,\ldots,0,1,0,\ldots,0)$ is nothing else than the principle of
a single big jump in the theory of subexponential distributions,
see Foss et al. (2011, Section 3.1).

An equivalent way is to consider the $L_\alpha$-norm
$g(\etaa)=(|\eta_1|^\alpha +\ldots+|\eta_d|^\alpha)^{1/\alpha}$
of Gaussian vector $\etaa$ which delivers an example of 
Gaussian chaos of order $1$.
It may be naturally extended for a general Minkowski functional
$h:\R^d\to\R^+$ where $h(\etaa)$ is again a Gaussian chaos of order $1$.
Earlier tail behavior of Minkowski's type of Gaussian chaos 
was studied by Pap and Richter (1988).

\textbf{Example 2.}
(Product of two Gaussian random variables $\xi_1$ and $\xi_2$)
Here we consider the case $d=2$ and assume without loss of
generality that $\mathrm{Var}\xi_1=\mathrm{Var}\xi_2=1$.
Denote the correlation coefficient by $\rho$, $\rho\neq-1$. Then
\[
B=\biggl(
\begin{array}
[c]{cc}
1 & \rho\\
\rho & 1
\end{array}
\biggr),\ \ \sqrt B=\frac{\sqrt{1+\rho}}{2}\biggl(
\begin{array}
[c]{cc}
1 & 1\\
1 & 1
\end{array}
\biggr) +\frac{\sqrt{1-\rho}}{2}\biggl(
\begin{array}
[c]{cc}
1 & -1\\
-1 & 1
\end{array}
\biggr),
\]
and $\xii$ has the same distribution as $\sqrt B\etaa$.
We consider the product
\[
h(\xi_1,\xi_2):=\xi_1\xi_2=g(\eta_1,\eta_2) =\frac{\rho}{2}
(\eta_1^2+\eta_2^2)+\eta_1\eta_2,
\]
so that $\alpha=2$ and  $g(\boldsymbol{u})=\rho(u_1^2+u_2^2)/2+u_1u_2$. Given $u_1^2+u_2^2=1$, the maximum of
$u_1u_2$ is attained on
$\mathcal M=\{(1/\sqrt 2,1/\sqrt 2),(-1/\sqrt 2,-1/\sqrt 2)\}$
and equals $\hat g=(1+\rho)/2$. At both points of the
maximum we have $g''_{2-1}(\boldsymbol{v}_j)=\rho-1$.
Calculating $h_0$ we obtain for $\rho\neq-1$, as $x\to\infty$,
\begin{align*}
\Pr\{\xi_1\xi_2>x\} & = \frac{1+\rho}{\sqrt{2\pi}}
x^{-1/2}e^{-x/(1+\rho)}\bigl(1+ O(1/x)\bigr),\\
p_{\xi_1\xi_2}(x) & = \frac{1}{\sqrt{2\pi}} x^{-1/2}
e^{-x/(1+\rho)}\bigl(1+O(1/x)\bigr).
\end{align*}

\textbf{Example 3.}
(Product of independent Gaussian random variables)
Let $\etaa=(\eta_1,\ldots,\eta_d)$ be a standard Gaussian standard vector,
that is, its components are $N(0,1)$ independent random variables.
Taking $g(\boldsymbol{u})=u_1\ldots u_d$ we  have $\alpha=d$ 
and further $\hat g=1/d^{d/2}$ since
\begin{align*}
\mathcal M =\{(\pm1/\sqrt d,\ldots,\pm1/\sqrt d)
\mbox{ with even number of negative coordinates}\},
\end{align*}
which  consists of $2^{d-1}$ points (the product $u_1\ldots u_d$
should be positive). Further, in
the spherical coordinates
\begin{align*}
g(\phii) =& \sin^{d-1}\varphi_1\ldots
\sin\varphi_{d-1}\cos\varphi_1\ldots\cos\varphi_{d-1}.
\end{align*}
For instance, at the point $(1/\sqrt d,\ldots,1/\sqrt d)$
we have $\cos\varphi_i=\sqrt{\frac{1}{d-i+1}}$ and
$\sin\varphi_i=\sqrt{\frac{d-i}{d-i+1}}$,
so that $\det J(1,\phii)=\sqrt{(d-1)!}/d^{(d-2)/2}$ at these points.
Additional calculations show,
at any point $\phii\in\mathcal M_\varphi$,
\[
g_{\varphi_i\varphi_i}''(\phii)=-2g(\phii)(d-i+1)
=-\frac{2(d-i+1)}{d^{d/2}},
\quad g_{\varphi_i\varphi_j}''(\phii)=0
\quad \text{for }i\neq j,
\]
which yields $|\det g''(\phii)|=2^{d-1}
d!/d^{d(d-1)/2}$. In this way we get the following answer
\[
p_{\eta_1\ldots\eta_d}(x)=\frac{2^{(d-1)/2}}
{\sqrt{2\pi d}}x^{1/d-1}e^{-dx^{2/d}/2}\bigl(1+O(x^{-2/d})\bigr)
\quad\mbox{ as }x\to\infty.
\]
The intuition behind this asymptotic behaviour is the following
(see e.g., Sornette (1998)): asymptotically, the tail of the product is
controlled by the realisations where all terms are of the same
order; therefore $p_{\eta_1\ldots\eta_d}(x)$ is,
up to the leading order, just the product of the $d$ marginal density functions,
evaluated at $x^{1/d}$.

\textbf{Example 4.}
(Product of components of a Gaussian vector) Let
$\xii=(\xi_1,\ldots,\xi_d)=\sqrt B\etaa$ be a Gaussian vector with
mean zero and with covariance matrix $B$ and consider
$h(\boldsymbol{u})=u_1\ldots u_d$.  Further, decompose the
symmetric positive-semidefinite matrix $\sqrt B$ as
$\sqrt B=Q^{\mathrm{T}}DQ$ where $Q$ is an orthogonal matrix
(the rows of which are eigenvectors of $\sqrt B$), and $D$ is
diagonal (having the eigenvalues of $\sqrt B$ on the diagonal).
Making use of the representation \eqref{eta.polar} we deduce
\[
h(\xii)=h(\sqrt B\etaa) =\chi^\alpha h(\sqrt B\zetaa)
=\chi^\alpha h(Q^{\mathrm{T}}DQ\zetaa).
\]
Since $Q$ is orthogonal the random vector
$\zetaa^*:=Q\zetaa$ is uniformly distributed on the
unit sphere $\Sp_{d-1}$. Therefore, $D\zetaa^*$ is
distributed on the ellipsoid $E$ with the semi-principal axes
of lengths equal to the diagonal elements of $D$.
The product of coordinates of $Q^{\mathrm{T}}\boldsymbol v$ on
$\boldsymbol v\in E$ has only finite number of points of maximum;
as above, denote this maximum by $\hat g$. It is not clear how to
identify the set $\mathcal M$ and the constants $\hat g$ and
$h_0$ explicitly, in terms of the covariance
matrix $B$, but we may guarantee that,
by Theorem \ref{th:m=0} with $r=1$
\[
p_{\xi_1\ldots\xi_d}(x) =\mathrm{const}\cdot x^{1/d-1}
e^{-(x/\hat g)^{2/d}/2} \bigl(1+O(x^{-2/d})\bigr)
\ \mbox{ as }x\to\infty.
\]

\textbf{Example 5.}
(Quadratic forms of independent $N(0,1)$ random variables)
Let $\etaa=(\eta_1,\ldots,\eta_d)$ be as in the previous section
and let $g(\etaa)=\sum_{i=1}^d a_i\eta_i^2$
where the constants $a_i\in\R$ are such that
\[
a_1\le a_2\le\ldots\le a_{d-m}<a_{d-m+1}=\ldots=a_d=a,\quad a>0.
\]
If $m=d$, then $g(\etaa)/a$ is a chi-square random variable.
So we consider the case $m\le d-1$. Since
\[
g(\boldsymbol{u}) = \sum_{i=1}^{d-m}a_i u_i^2+ a\sum_{i=d-m+1}^d u_i^2
\]
and all $a_i<a$ for $i\le d-m$, the maximum of $g(\boldsymbol{u})$
given $\|\boldsymbol{u}\|=1$
is attained at any point $\boldsymbol{u}$
such that $u_{d-m+1}^2+\ldots+u_d^2=1$ and $u_1=\ldots=u_{d-m}=0$ implying  $\hat g=a$. We have further
\begin{align*}
g(\phii) & = a_1\cos^2\varphi_1
+\sum_{i=2}^{d-m}a_i\sin^2\varphi_1\ldots
\sin^2\varphi_{i-1}\cos^2\varphi_i
+a\sin^2\varphi_1\ldots\sin^2\varphi_{d-m}.
\end{align*}
The set $\mathcal M_\varphi$ of dimension $m-1$
is the sub-parallelepiped of $\Pi_{d-1}$, namely,
${\mathcal M}_\varphi=\{\phii\in\Pi_{d-1}:\varphi_1=\ldots=\varphi_{d-m}=\pi/2\}$
for $m\ge 2$, its inverse image $\mathcal M$ is a unit sphere
$\Sp_{m-1}$; in the case $m=1$ it consists of two points
$(\pi/2,\ldots,\pi/2,\pi/2)$ and $(\pi/2,\ldots,\pi/2,3\pi/2)$.
For $\phii\in\mathcal M_\varphi$, the matrix
$g_{d-1}''(\phii)$ is diagonal with first entries
$2(a_i-a)$ for $i=1$, \ldots, $d-m$
and zeros on the rest of diagonal, so
\[
|\det g''_{d-m}(\phii)|=2^{d-m}\prod_{i=1}^{d-m}(a-a_i)
\]
does not depend on $\phii\in\mathcal M_\varphi$. Therefore,
$$
\int_{\mathcal M_\varphi}
\frac{|\det J(1,\phii)|}{\sqrt{|\det g''_{d-m}(\phii)|}}{\rm d}\Vfi
=\frac{{\rm mes}\,\mathcal M}{2^{\frac{d-m}{2}}
\prod_{i=1}^{d-m}\sqrt{a-a_i}}.
$$
Taking into account that
${\rm mes}\,\mathcal M={\rm mes}\,\Sp_{m-1}=2\pi^{m/2}/\Gamma(m/2)$
we finally deduce, as $x\to\infty$
\begin{align}\label{Z}
p_{\sum_{i=1}^d a_i \eta_i^2}(x)
& = \frac{1}{a2^{m/2}\Gamma(m/2)}
\prod_{i=1}^{d-m}\frac{1}{\sqrt{1-a_i/a}} (x/a)^{m/2-1}
e^{-x/2a}(1+O(1/x)),
\end{align}
which agrees (for the first order asymptotics) with Hoeffding ({1964}) 
(see also Zolotarev (1961), Imkeller (1994), Piterbarg (1994, 1996),
H\"usler et al. (2002)). 

\textbf{Example 6.}
(Scalar product of Gaussian random vectors) Closely
related to Example 5 is the scalar product of two independent
Gaussian random vectors, namely we consider the Gaussian chaos
$g(\etaa, \boldsymbol{\eta}^*)=\sum_{i=1}^d a_i\eta_i\eta_i^*$
with $\eta_i$, $\eta_i^*$, $i\le d$, independent $N(0,1)$ random
variables. Indeed, since $\eta_i\eta_i^*$ coincides in distribution with
\[
\frac{\eta_i+\eta_i^*}{\sqrt2}\frac{\eta_i-\eta_i^*}{\sqrt2}
=\frac{\eta_i^2-\eta_i^{*2}}{2}
\]
we have the equality in distribution
\[
g(\etaa,\etaa^*) \overset{d}{=}\frac{1}{2}\biggl(\sum
_{i=1}^d a_i\eta_i^2 -\sum_{i=1}^d a_i\eta_i^{*2}\biggr).
\]
Therefore, if
$$
|a_1|\le |a_2|\le\ldots\le |a_{d-m}|<|a_{d-m+1}|
=\ldots=|a_d|=a,\quad a>0,
$$
then the asymptotics of the density given by \eqref{Z} is applicable,
and we have as $x\to\infty$
\begin{align*}
p_{g(\etaa,\etaa^*)}(x)
& = \frac{1}{a2^{d/2}\Gamma(m/2)}
\prod_{i=1}^{d-m}\frac{1}{\sqrt{1-a_i^2/a^2}}
(x/a)^{m/2-1} e^{-x/a}(1+O(1/x)).
\end{align*}
We note that results for the scalar products
of Gaussian random variables are derived in
Ivanoff and Weber (1998) and Hashorva et al. (2012).

\textbf{Example 7.}
(Determinant of a random Gaussian matrix)
Let $A=[A_{ij}]_{i,j=1}^n$ be a random square matrix of
order $n$ whose entries $A_{ij}$ are independent $N(0,1)$
random variables. Then its determinant is the following
function of $A_{ij}$:
$$
\det A=\sum_{\sigma\in S_n}{\rm sgn}(\sigma)\prod_{i=1}^n A_{i,\sigma_i},
$$
where $S_n$ is the set of all permutations of the set
$\{1,2,\ldots,n\}$ and ${\rm sgn}(\sigma)$ denotes the signature
of $\sigma\in S_n$. Clearly, the determinant $g(A):=\det A$
is a continuous homogeneous function of order $\alpha=n$.
Here we have $d=n^2$.

The determinant of the matrix $A$ represents the (oriented)
volume of the parallelepiped generated by the vectors
${\boldsymbol A}_i:=(A_{i1},\ldots,A_{in})$, $i=1$, \ldots, $n$.
Given
$$
\sum_{i=1}^n\|{\boldsymbol A}_i\|^2=1,
$$
the maximal volume of this parallelepiped is attained on
orthogonal vectors ${\boldsymbol A}_i$ which are of the same length,
that is, on the $n$-dimensional cube with side of length $1/\sqrt n$.
Therefore,
$$
\hat g:=\max_{A:\sum_{i,j=1}^n A_{ij}^2=1} \det A=n^{-n/2}.
$$
The manifold consisting of points
where the maximum $\hat g$ of $g(A)$, $A\in\Sp_{n^2-1}$,
is attained, that is,
$$
\mathcal M:=\bigl\{A:\det A>0,\|{\boldsymbol A}_1\|=
\ldots=\|{\boldsymbol A}_n\|=1/\sqrt n\mbox{ and }
{\boldsymbol A}_1,\ldots,{\boldsymbol A}_n
\mbox{ are orthogonal}\bigr\}
$$
has dimension $m=(n^2-n)/2$. Therefore, by Theorem \ref{th:m=1.d-2},
\begin{eqnarray} \label{eqdet0}
\Pr\{\det A>x\} &=&
cx^{\frac{n-1}{2}-\frac{1}{n}}e^{-nx^{2/n}/2}(1+O(x^{-2/n}))
\quad\mbox{ as }x\to\infty,
\end{eqnarray}
for some $c=c(n)>0$; the computation of this constant is questionable.
This answer agrees with Theorem 10.1.4(i) by Barbe (2003)
in the exponential term and gives the correct power term.

Another way to show this result is to recall from Pr\'ekopa (1967, Theorem 2) that
\begin{equation}\label{eqdet} 
\det A\stackrel{d}{=} \prod_{i=1}^n \chi_i^2, 
\end{equation}
where $\chi_1^2$, \ldots, $\chi_n^2$ are independent random variables 
and $\chi_i^2$ is chi-square distributed with $i$ degrees of freedom.
In the case $n=2$ it easily follows by conditioning on $\eta_1$ and $\eta_3$:
$$ 
\det A^2 \stackrel{d}{=} (\eta_1\eta_2+\eta_3\eta_4)^2\stackrel{d}{=} 
\eta_5^2(\eta_6^2+\eta_7^2)=: \chi_1^2 \chi_2^2,
$$
where $\eta_1$, \ldots, $\eta_7$ are independent $N(0,1)$ random variables. 
The representation (\ref{eqdet}) provides an alternative way of 
deducing the tail asymptotics of $\det A$, 
since we can readily apply Lemma 3.2 in \cite{Are}.

\textbf{Example 8.}
(Gaussian orthogonal ensemble) 
Now let $A=[A_{ij}]_{i,j=1}^n$ be a random square
symmetric matrix of order $n$ whose random entries $A_{ij}$
are independent for $1\le i\le j\le n$.
Let $A_{ij}=\eta_{ij}$ for $j>i$ and $A_{ii}=\sigma\eta_{ii}$
where $\eta_{ij}$ are independent standard random variables and $\sigma>1$.
In the special case $\sigma=\sqrt 2$ the matrix $A$ is called the 
{\it Gaussian orthogonal ensemble}.

Here the determinant $g(A):=\det A$ is again a continuous
homogeneous function of order $\alpha=n$; $d=(n^2+n)/2$.
Due to the coefficients $\sigma>1$ on the diagonal,
the maximal volume of the corresponding parallelepiped is
attained on the orthogonal vectors $(\pm\sigma/\sqrt{n},0,\ldots,0)$, \ldots,
$(0,\ldots,0,\pm\sigma/\sqrt{n})$ with even number of minuses.
Hence, $\hat g=(\sigma^2/n)^{n/2}$. Since $\mathcal M$ is finite,
we apply Theorem \ref{th:m=0} and deduce that
\begin{eqnarray*}
\Pr\{\det A>x\} &\sim& cx^{-1/n}e^{-nx^{2/n}/2\sigma^2}
\end{eqnarray*}
as $x\to\infty$, for some $c=c(n)>0$.

An alternative approach for computing asymptotics of the
tail of the Gaussian orthogonal ensemble (where $\sigma=\sqrt 2$) 
is to make use of the fact that the joint density function of the eigenvalues
$\lambda_1(A)\le\ldots\le\lambda_n(A)$ is known and is equal to
$$
c'e^{-\|y\|^2/4}{\mathbb I}\{y_1\le\ldots\le y_n\}\prod_{i<j}(y_j-y_i),
$$
with some explicitly known normalising constant $c'=c'(n)>0$,
see e.g., Theorem 2.5.2 in Anderson et al. (2010).
Clearly this approach is more complicated from computational
point of view because of the singularity of the product
$\prod_{i<j}(y_j-y_i)$ on the diagonal $y_1=\ldots=y_n$.

Indeed, there are Gaussian chaoses where Theorems \ref{th:m=0} 
and \ref{th:m=1.d-2} are not straightforward applicable because of 
degeneracy of their Hessian on the set $\mathcal M$ of extremal points. 
This is exactly the case of diameter of a random Gaussian chaos 
which is discussed next.

\textbf{Example 9.}
(The diameter of a random Gaussian cloud) 
Let $\etaa_k=(\eta_{k1},\ldots,\eta_{km})$, $k=1$, \ldots, $n$, 
be i.i.d. random vectors in $\R^m$; here $\eta_{kl}$, $k=1$, \ldots, $n$,
$l=1$, \ldots, $m$, are independent $N(0,1)$ random variables.
The set of random points $\{\etaa_k,k\le n\}$ 
may be called the {\it Gaussian cloud}. 
The problem is how to approximate the distribution of its {\it diameter}
\begin{eqnarray*}
D_n &=& \max_{1\le k\le l\le n}\|\etaa_k-\etaa_l\|.
\end{eqnarray*}
In \cite{MR}, Matthews and Rukhin study limit behavior of $D_n^2$ as $n\to\infty$.
Here we discuss the problem of estimation of the tail of $D_n$ for a fixed $n$.
First of all notice that it is equivalent to tail estimation of 
$g(\etaa_1,\ldots,\etaa_n):=D_n^2$ which represents
a smooth Gaussian chaos of order $\alpha=2$, with $d=mn$.
Since the cases $m=1$ and $m\ge 2$ are different, consider them separately.

First consider the case of dimension $1$, $m=1$.
For any $k\not=l$, introduce $T_{kl}$ as the set of all vectors 
$(v_1,\ldots,v_n)\in\R^n$ such that $v_i=0$ for all $i\not\in\{k,l\}$, 
$v_k=-v_l$ and $v_k=\pm1/\sqrt 2$. Then, given
$v_1^2+\ldots+\eta_n^2=1$,
the maximal value $\hat g$ of $D_n^2$ is attained on 
$\mathcal M=\cup_{k\not=l}T_{kl}$; this set consists of $n(n-1)$ points. 
In particular,
$$
\hat g=\max_{\sum_{k=1}^n v_k^2=1} \max_{1\le k\le l\le n} (v_k-v_l)^2
=(1/\sqrt 2-(-1/\sqrt 2))^2 = 2.
$$
Therefore, by Theorem \ref{th:m=0},
\begin{eqnarray*}
\Pr\{D_n^2>x\} &=& h_0(x/2)^{-1/2}e^{-x/4}(1+O(1/x))\quad\mbox{ as }x\to\infty,
\end{eqnarray*}
where
\begin{eqnarray*}
h_0 &:=& \frac{1}{\sqrt{2\pi}}
\sum_{(v_1,\ldots,v_n)\in\cup_{k,l} T_{kl}}
\Bigl|\det\Bigl(\frac{g''_{d-1}(v_1,\ldots,v_n)}{4}-I_{d-1}\Bigr)\Bigr|^{-1/2}.
\end{eqnarray*}
The latter sum consists of $n(n-1)$ equal terms. 
Consider a typical representative, 
$\boldsymbol V_0:=(1/\sqrt 2,-1/\sqrt 2,,0,\ldots,0)$, 
which contains $n-2$ zeros. Consider the following orthogonal
system of vectors $\boldsymbol E_1$, \ldots, $\boldsymbol E_{n-1}$ 
in the hyperplane 
$\mathcal L:=\{\boldsymbol V\in\R^n: (\boldsymbol V,\boldsymbol V_0)=0\}$:
the vector $\boldsymbol E_1:=(1/\sqrt 2,1/\sqrt 2,0,\ldots,0)\in\R^n$ 
plus an orthogonal system $\boldsymbol E_2$, \ldots, $\boldsymbol E_{n-1}$
in $\{(0,0,v_3,\ldots,v_n)\in\R^n\}$.
Since the function $g(v_1,\ldots,v_n)$ is equal to $(v_1-v_2)^2$
in some neighborhood of the point $\boldsymbol V_0$,
the Hessian of the function $g$ at point $\boldsymbol V_0$ is 
the following square matrix of size $n$
\begin{eqnarray*}
g'' &=& \left(
\begin{array}{rrccc}
2 & -2 & 0 &\ldots & 0\\
-2 & 2 & 0 &\ldots & 0\\
0 & 0 & 0 &\ldots & 0\\
 &  & \ldots & &\\
0 & 0 & 0 &\ldots & 0
\end{array}
\right).
\end{eqnarray*}
Then the Hessian of the function $g$ restricted to the hyperplane 
$\boldsymbol V_0+\mathcal L$ is zero square matrix of size $n-1$, 
because its entries are equal to $(g''\boldsymbol E_i,\boldsymbol E_j)$. 
Hence we conclude that $h_0=\frac{1}{\sqrt{2\pi}} n(n-1)$, 
so that in dimension $1$
\begin{eqnarray*}
\Pr\{D_n>x\} &=& \Pr\{D_n^2>x^2\}
= \frac{n(n-1)}{\sqrt\pi x} e^{-x^2/4}(1+O(1/x^2))
\quad\mbox{ as }x\to\infty.
\end{eqnarray*}

Next, we  show that in dimension greater than $1$ the situation is more complicated.
For any $k\not=l$, introduce $T_{kl}$ as the set of all vectors 
$(\boldsymbol v_1,\ldots,\boldsymbol v_n)\in\R^{mn}$ such that 
$\boldsymbol v_i=\boldsymbol 0\in\R^m$ for all $i\not\in\{k,l\}$, 
$\boldsymbol v_k=-\boldsymbol v_l\in\R^m$ and 
$\boldsymbol v_{ki}=\pm1/\sqrt{2m}$ for all $i\le m$. Then, given
$$
\|\boldsymbol v_1\|^2+\ldots+\|\boldsymbol v_n\|^2=1,
$$
the maximal value $\hat g$ of $D_n^2$ is attained on 
$\mathcal M=\cup_{k\not=l}T_{kl}$; this set consists of 
$2^m\frac{n(n-1)}{2}$ points. In particular,
$$
\hat g=\max_{\sum_{k=1}^n\sum_{i=1}^m \boldsymbol v_{ki}^2=1} 
\max_{1\le k\le l\le n} \|\boldsymbol v_k-\boldsymbol v_l\|^2
=\sum_{i=1}^m (1/\sqrt{2m}-(-1/\sqrt{2m}))^2 = 2
$$
independently of $m$ and $n$. In order to apply Theorem \ref{th:m=0},
we need to compute the following constant
\begin{eqnarray*}
h_0 &:=& \frac{1}{\sqrt{2\pi}}
\sum_{(\boldsymbol v_1,\ldots,\boldsymbol v_n)\in\cup_{k,l} T_{kl}}
\Bigl|\det\Bigl(\frac{g''_{d-1}(\boldsymbol v_1,\ldots,\boldsymbol v_n)}{4}-I_{d-1}\Bigr)\Bigr|^{-1/2}.
\end{eqnarray*}
The latter sum consists of equal terms. 
Consider a typical representative, 
$$
\boldsymbol V_0:=(\boldsymbol 1/\sqrt{2m},-\boldsymbol 1/\sqrt{2m},
\boldsymbol 0,\ldots,\boldsymbol 0),
$$ 
which contains $m$ coordinates equal to $1/\sqrt{2m}$, 
$m$ coordinates equal to $-1/\sqrt{2m}$, and $m(n-2)$ zeros. 

Consider the following orthogonal system of vectors 
$\boldsymbol E_1$, \ldots, $\boldsymbol E_{mn-1}$ in the hyperplane 
$\mathcal L:=\{\boldsymbol V\in\R^{mn}: (\boldsymbol V,\boldsymbol V_0)=0\}$:
\begin{eqnarray*}
\left(\begin{array}{l}
\boldsymbol E_1\\
\boldsymbol E_2\\
\ldots\\
\boldsymbol E_{m-1}\\
\boldsymbol E_m\\
\boldsymbol E_{m+1}\\
\ldots\\
\boldsymbol E_{2m-2}\\
\boldsymbol E_{2m-1}
\end{array}
\right) &=& \left(\begin{array}{lllll}
\boldsymbol e_1 & \boldsymbol 0 & \boldsymbol 0 & \ldots & \boldsymbol 0\\
\boldsymbol e_2 & \boldsymbol 0 & \boldsymbol 0 & \ldots & \boldsymbol 0\\
\ldots\\
\boldsymbol e_{m-1} & \boldsymbol 0 & \boldsymbol 0 & \ldots & \boldsymbol 0\\
\boldsymbol 0 & \boldsymbol e_1 & \boldsymbol 0 & \ldots & \boldsymbol 0\\
\boldsymbol 0 & \boldsymbol e_2 & \boldsymbol 0 & \ldots & \boldsymbol 0\\
\ldots\\
\boldsymbol 0 & \boldsymbol e_{m-1} & \boldsymbol 0 & \ldots & \boldsymbol 0\\
\boldsymbol 1/\sqrt{2m} & \boldsymbol 1/\sqrt{2m} & \boldsymbol 0 & \ldots & \boldsymbol 0
\end{array}
\right)
\end{eqnarray*}
where $\boldsymbol e_k=\frac{1}{\sqrt{k(k+1)}}(1,\ldots,1,-k,0,\ldots,0)$ with $k$ units and 
$m-k-1$ zeros; plus an orthogonal system 
$\boldsymbol E_{2m}$, \ldots, $\boldsymbol E_{mn-1}$
in $\{(\boldsymbol 0,\boldsymbol 0,v_{2m+1},\ldots,v_{mn})\in\R^{mn}\}$. Since 
$g(\boldsymbol v_1,\ldots,\boldsymbol v_n)=\|\boldsymbol v_1-\boldsymbol v_2\|^2$
in some neighborhood of the point $\boldsymbol V_0$,
the Hessian of the function $g$ at point $\boldsymbol V_0$ is 
the following square matrix of size $mn$
\begin{eqnarray*}
g'' &=& 2\left(
\begin{array}{rrccc}
I_m & -I_m & & & \\
-I_m & I_m & & & \\
 & & 0 & & \\
 & & & \ddots &\\
 & & & & 0
\end{array}
\right),
\end{eqnarray*}
with $m(n-2)$ zero diagonal entries.
Computing the entries of the Hessian of the function $g$ restricted 
to the hyperplane $\boldsymbol V_0+\mathcal L$ via 
$(g''\boldsymbol E_i,\boldsymbol E_j)$ we get that 
\begin{eqnarray*}
g''_{d-1} = 2\left(
\begin{array}{rrlll}
I_{m-1} & -I_{m-1} & & & \\
-I_{m-1} & I_{m-1} & & & \\
 & & 0 & &\\
 & & & \ddots &\\
 & & & & 0
\end{array}
\right),
&&
\frac{g''_{d-1}}{4}-I_{d-1} = \frac{1}{2}\left(
\begin{array}{rrlll}
-I_{m-1} & -I_{m-1} & & & \\
-I_{m-1} & -I_{m-1} & & & \\
 & & -2 & &\\
 & & & \ddots &\\
 & & & & -2
\end{array}
\right).
\end{eqnarray*}
Hence we conclude that $\det(g''_{d-1}/4-I_{d-1})=0$ in the case $m\ge 2$,
so that $\Pr\{D_n>x\}x e^{x^2/4}\to\infty$ as $x\to\infty$.
This example calls for study of degenerated Hessians but we do not
concern this question in the current paper.

Similar examples can be given for spherical  chaos by applying our Theorem 3. 
The calculation of $h_0$ therein is readily obtained using
the results of previous examples. 
Note that the determination of $\hat g$ and the parameter $m$ is the same 
as for the Gaussian chaos. In order to avoid repetition we present 
only the case of Examples 1 and 7; the remaining cases 
can be easily extended by studying the next two examples. 

\textbf{Example 10.} Let  $\etaa$ be $d$-dimensional random vector which 
is spherically distributed, such that (1) holds with $\chi$ 
a positive random radius. Consider 
$g(\etaa)=|\eta_1|^\alpha +\ldots+|\eta_d|^\alpha $ with $\alpha>0$ given 
and assume that $\chi^\alpha$ has distribution function in the Gumbel 
max-domain of attraction with some scaling function $w$.

If $\alpha\in(0,2)$, then by Example 1 we have $m=0$ and 
$\hat g=d^{1-\alpha/2}$, hence by Theorem \ref{thPolar} we find that, 
as $x\to \infty$,
\begin{eqnarray*}
\Pr\{g(\etaa)> x\} &\sim& 
h_0(xw(x/\hat g))^{\frac{1-d}{2}}\Pr\{\chi^\alpha>x/\hat g\}
\end{eqnarray*}
as $x\uparrow \hat g x_+$, where
\begin{eqnarray*}
h_0=\frac{2^{3d/2-3/2}\Gamma(d/2)\hat g^{\frac{d-1}{2}}}
{\sqrt\pi(\alpha(2-\alpha))^{\frac{d-1}{2}}}.
\end{eqnarray*}

In the case $\alpha>2$, then $\hat g=1$ and $m=0$ as in Example 1
and we find that 
\begin{eqnarray*}
\Pr\{g(\etaa)> x\} &\sim& 
\Bigl(\frac{2}{\alpha}\Bigr)^{\frac{d-1}{2}}\frac{d\Gamma(d/2)}{\sqrt\pi}
\bigl(xw(x)\bigr)^{\frac{1-d}{2}}\Pr\{\chi^\alpha>x\}
\quad\mbox{as }x\to\infty.
\end{eqnarray*} 
We note in passing that the above results agree with the direct 
calculations in \cite{HE4}; the case $\alpha=2$ is discussed in \cite{HE0}.

\textbf{Example 11.} (Determinant of a random spherical  matrix)
Let $A=[A_{ij}]_{i,j=1}^n$ be a random square matrix of
order $n$ and let $\boldsymbol A^*=(\boldsymbol a_1,\ldots,\boldsymbol a_n)$ 
be the $n^2\times 1$ vector obtained pasting the rows of $A$, i.e., 
$\boldsymbol a_i=(A_{i1},\ldots,A_{in})$ is the $i$th row of $A$. 
Suppose that $\boldsymbol A^*$ is a spherically 
symmetric random vector meaning that 
$$
\boldsymbol A^*\stackrel{d}{=} \chi\zetaa
$$
with $\chi>0$ being independent of $\zetaa$ which is uniformly distributed 
on $\Sp_{n^2-1}$. We consider again $g(A)=\det A$ which is a continuous 
homogeneous function of order $\alpha=n$. 
Note that if $\chi^2$ has a chi-square distribution with $n^2$ degrees 
of freedom, then $A$ is the matrix in Example 7. Hence for this case, 
if $\chi^\alpha \in GMDA(w,x_+)$ with $x_+=\infty$, say,
then since by Example 7 we have $m=(n^2-n)/2$ and $\hat g=n^{-n/2}$, 
$d=n^2$, then Theorem 3 entails
\begin{eqnarray}
\Pr\{\det A>x\} &\sim& 
c^*\bigl(xw(xn^{n/2})\bigr)^{-\frac{n^2+n-2}{4}}\Pr\{\chi^n>xn^{n/2}\}
\end{eqnarray}
as $x\to\infty$ for some constant $c^*$, which can be calculated iteratively 
by applying Lemma 3.2 in \cite{Are} as mentioned in Example 7.

\section{Proof of Theorem \ref{th:m=0}}
\label{sec:proof.0}

The proof is based on the polar representation \eqref{eta.polar}
for a $d$-dimensional centered Gaussian random vector $\etaa$
with identity covariance matrix, $\etaa\overset{d}=\chi\zetaa$,
where $\chi$ and $\zetaa$ are independent,
$\chi^2=\sum_{i=1}^d\eta_i^2$ has $\chi^2$-distribution with 
$d$ degrees of freedom and $\zetaa$ is uniformly distributed 
on the unit sphere $\Sp_{d-1}\subset\R^d$. 
The tail distribution of the random variable
$$
g(\etaa)\stackrel{d}{=} g(\chi\zetaa)
=\chi^\alpha g(\zetaa)
$$
is equal to
\begin{eqnarray}\label{multipl.our}
\Pr\{g(\etaa)>x\}
&=& \int_{x/\hat g}^\infty p_{\chi^\alpha}(y)\Pr\{g(\zetaa)>x/y\} dy\nonumber\\
&=& \frac{1}{\alpha2^{d/2-1}\Gamma(d/2)}
\int_{x/\hat g}^\infty y^{d/\alpha-1}e^{-y^{2/\alpha}/2}
\Pr\{g(\zetaa)>x/y\}{\rm d}y
\end{eqnarray}
by the equality \eqref{chi.alpha} and boundedness $g(\zetaa)\le\hat g$.
In order to compute the asymptotics for the latter integral,
we first need to estimate the probability
$\Pr\{g(\zetaa)>\hat g-t\}$ for small positive values of $t$. Hereinafter ${\rm Vol}\,\mathbb B_{d-1}$ stands for the volume
of the unit ball $\mathbb B_{d-1}$ in $\R^{d-1}$.

\begin{lemma}\label{asy.for.g.mu.0}
Under the conditions of Theorem \ref{th:m=0}
$$
\Pr\{g(\zetaa)>\hat g-t\} =
\sum_{i=0}^r g_i t^{\frac{d-1}{2}+i}
+o(t^{\frac{d-1}{2}+r})\quad\mbox{as }t\downarrow 0,
$$
where
$$
g_0=2^{\frac{d-1}{2}}\frac{{\rm Vol}\,\mathbb B_{d-1}}{{\rm mes}\,\Sp_{d-1}}
\sum_{j=1}^k\bigl|\det\bigl(g_{d-1}''(\boldsymbol v_j)
-(\alpha\hat g)I_{d-1}\bigr)\bigr|^{-1/2},
$$
and where further coefficients $g_1$, \ldots, $g_r$ only depend on $\alpha$,
$\hat g$, and derivatives of $g(\phii)$ at points $\phii_j$.
\end{lemma}

Since ${\rm Vol}\,\mathbb B_{d-1}=\frac{\pi^{(d-1)/2}}{\Gamma((d+1)/2)}$
and ${\rm mes}\,\mathbb S_{d-1}=\frac{2\pi^{d/2}}{\Gamma(d/2)}$, 
the expression for the coefficient $g_0$ may be rewritten
in the hyperspherical coordinates as follows (see \eqref{cart.local}):
$$
g_0=\frac{2^{d/2-1}}{\sqrt{2\pi}}
\frac{\Gamma(d/2)}{\Gamma((d+1)/2)}
\sum_{j=1}^k\frac{|\det J(1,\phii_j)|}{\sqrt{|\det g''(\phii_j)|}}.
$$

\proof
Without loss of generality we consider the case where $\mathcal M$
consists of a single point $\boldsymbol v_1$. First prove that
$$
\Pr\{g(\zetaa)>\hat g-t\} \sim g_0 t^{\frac{d-1}{2}}
\quad\mbox{as }t\downarrow 0.
$$

Introduce the hyperspherical coordinates of the random vector
$\zetaa=(\zeta_1,\ldots,\zeta_d)\in\Sp_{d-1}$ as
$\nuu=(\nu_1,\ldots,\nu_{d-1})\in\Pi_{d-1}$.
Since $\zetaa$ is uniformly distributed on the unit sphere $\Sp_{d-1}$
in $\R^d$, the density function of the random vector
$\nuu=(\nu_1,\ldots,\nu_{d-1})\in\Pi_{d-1}$
equals $\frac{|\det J(1,\phii)|}{{\rm mes}\,\Sp_{d-1}}$,
$\phii\in\Pi_{d-1}$, which implies
\begin{eqnarray}\label{deco.repr}
\Pr\{g(\zetaa)>\hat g-t\} &=& \frac{1}{{\rm mes}\,\Sp_{d-1}}
\int_{\phii\in\Pi_{d-1}:g(\phii)>\hat g-t} |\det J(1,\phii)| d\phii.
\end{eqnarray}
Since $g$ is at least twice differentiable and attains its maximum
at point $\phii_1$, 
\begin{equation*}
g(\phii)=\hat g+\frac12 \Bigl((g''(\phii_1)+A(\phii))(\phii-\phii_1),\phii-\phii_1 \Bigr),
\end{equation*}
where all the coefficients of the matrix $A(\phii)$
go to $0$ as $\phii\to\phii_1$.
Therefore, the inequality $g(\phii)>\hat g-t$ is equivalent to
$$
-\Bigl((g''(\phii_1)+A(\phii))(\phii-\phii_1),\phii-\phii_1 \Bigr) \le 2t.
$$
Fix $\varepsilon>0$. There exists $\delta>0$ such that 
\begin{equation}\label{A.I.eps}
-\varepsilon I_{d-1}\le A(\phii)\le \varepsilon I_{d-1}
\quad\mbox{for all } \phii\mbox{ such that } \|\phii-\phii_1\|\le\delta.
\end{equation}
Then, for all sufficiently small $t>0$, the set 
$\{\phii:g(\phii)>\hat g-t\}$ is contained in the
$(d-1)$-dimensional ellipsoid
$$
\Bigl((-g''(\phii_1)-\varepsilon I_{d-1})(\phii-\phii_1),\phii-\phii_1 \Bigr) \le 2t,
$$
whose volume is 
$$
\frac{{\rm Vol}\,\mathbb B_{d-1}}
{\sqrt{|\det(g''(\phii_1)+\varepsilon I_{d-1})|}}(2t)^{\frac{d-1}{2}}.
$$
On the other hand, \eqref{A.I.eps} implies that,
for all sufficiently small $t>0$, the set 
$\{\phii:g(\phii)>\hat g-t\}$ contains the
$(d-1)$-dimensional ellipsoid
$$
\Bigl((-g''(\phii_1)+\varepsilon I_{d-1})(\phii-\phii_1),\phii-\phii_1 \Bigr) \le 2t,
$$
whose volume is 
$$
\frac{{\rm Vol}\,\mathbb B_{d-1}}
{\sqrt{|\det(g''(\phii_1)-\varepsilon I_{d-1})|}}(2t)^{\frac{d-1}{2}}.
$$
Since $\varepsilon>0$ can be chosen arbitrary small,
the above arguments yield that
the volume of the set $\{\phii:g(\phii)>\hat g-t\}$ is proportional to
$$
\frac{{\rm Vol}\,\mathbb B_{d-1}}{\sqrt{|\det g''(\phii_1)|}}
(2t)^{\frac{d-1}{2}}\quad\mbox{ as }t\downarrow 0.
$$
Together with \eqref{deco.repr} this proves the required
asymptotic behavior of the probability $\Pr\{g(\zetaa)>\hat g-t\}$.

Next, given that $g$ is differentiable sufficiently many times,
the probability $\Pr\{g(\zetaa)>\hat g-t\}$ clearly possesses the
decomposition with terms $t^{\frac{d-1+i}{2}}$. It turns out that
in reality all terms with $i$ odd have zero coefficients. So,
it remains to prove that the asymptotic expansion of the integral
\eqref{deco.repr} only contains the terms $t^{\frac{d-1}{2}+i}$ 
and does not contain terms of order $\frac{d-1}{2}+\frac{1}{2}+i$. 
It is done in Lemmas \ref{l:1.f-1} and \ref{l:1.f-d} below and
the proof of Lemma \ref{asy.for.g.mu.0} follows.

\begin{lemma}\label{l:1.f-1}
Let a function $g(u):[-1,1]\to\R^+$ possess an asymptotic expansion
$$
g(u)=\sum_{i=2}^{2r+2} g_i u^i+o(u^{2r+2}) \quad\mbox{as }u\to 0,
$$
where $g_2>0$. Let $g(u)$ be strictly decreasing for 
$u\in[-1,0]$ and strictly increasing for $u\in[0,1]$.
For $t\in(0,g(-1)\wedge g(1))$, denote by $u^+(t)$ the unique 
positive value of $g^{-1}(t)$ and by $u^-(t)$ the negative one. 
Let a function $w(u)$ possess an asymptotic expansion
$$
w(u)=\sum_{i=0}^{2r} w_i u^i+o(u^{2r}) \quad\mbox{as }u\to 0.
$$
Then
\begin{eqnarray}\label{u.pm}
\int_{u^-(t)}^{u^+(t)}w(u')du' &=& \frac{2w_0}{\sqrt{g_2}}\sqrt t
+\sum_{i=1}^r u_i t^{1/2+i}+o(t^{1/2+r})
\quad\mbox{as }t\downarrow 0
\end{eqnarray}
and
\begin{eqnarray}\label{u.pm.abs}
\int_{u^-(t)}^{u^+(t)} |u'|w(u')du' &=& \frac{w_0}{g_2}t
+\sum_{i=2}^{r+1} \widetilde u_i t^i+o(t^{r+1})
\quad\mbox{as }t\downarrow 0,
\end{eqnarray}
where coefficients $u_1$, \ldots, $u_r$ and $\widetilde u_2$, \ldots,
$\widetilde u_{r+1}$ only depend on $g_i$'s and $w_i$'s.

Moreover, let $\Theta$ be a parameter set and let, 
for every fixed $\theta\in\Theta$, the functions
$g(u,\theta)$ and $w(u,\theta)$ satisfy the conditions stated above.
Suppose that 
\begin{eqnarray*}
\inf_{\theta\in\Theta}g_2(\theta) &>& 0
\end{eqnarray*}
and that all coefficients are uniformly bounded on $\Theta$,
\begin{eqnarray*}
\sup_{\theta\in\Theta}|g_i(\theta)| &<& \infty,\quad i\in\{1,2,\ldots,2r+2\},\\
\sup_{\theta\in\Theta}|w_i(\theta)| &<& \infty,\quad i\in\{0,1,\ldots,2r\},
\end{eqnarray*}
and the remainder terms are uniform on $\Theta$:
\begin{eqnarray*}
\sup_{\theta\in\Theta}\biggl|g(u,\theta)-\sum_{i=2}^{2r+2} 
g_i(\theta) u^i\biggr| &=& o(u^{2r+2}),\\
\sup_{\theta\in\Theta}\biggl|w(u,\theta)-\sum_{i=0}^{2r} 
w_i(\theta) u^i\biggr| &=& o(u^{2r})
\end{eqnarray*}
as $u\to 0$. Then, for 
$$
t\in\Bigl(0,\min_{\theta\in\Theta}g(-1,\theta)\wedge 
\min_{\theta\in\Theta} g(1,\theta)\Bigr),
$$ 
the asymptotic expansion
\begin{eqnarray}\label{u.pm.uni}
\int_{u^-(t,\theta)}^{u^+(t,\theta)} w(u',\theta)du' &=& 
\frac{2w_0(\theta)}{\sqrt{g_2(\theta)}}\sqrt t+
\sum_{i=1}^r u_i(\theta) t^{1/2+i}+o(t^{1/2+r})
\end{eqnarray}
holds as $t\downarrow 0$ uniformly on $\Theta$ and
\begin{eqnarray}\label{u.pm.abs.uni}
\int_{u^-(t,\theta)}^{u^+(t,\theta)} |u'|w(u',\theta)du' &=& 
\frac{w_0(\theta)}{g_2(\theta)} t+
\sum_{i=2}^{r+1} \widetilde u_i(\theta) t^i+o(t^{r+1}).
\end{eqnarray}

\end{lemma}

Conditions of Lemma \ref{l:1.f-1} almost immediately imply that
$$
u^\pm(t)=\sum_{i=0}^{2r} u^\pm_i t^{1/2+i/2}+o(t^{1/2+r})
\quad\mbox{as }t\downarrow 0,
$$
so that the asymptotic expansion \eqref{u.pm}, with $w(u)=u$, 
is equivalent to the nontrivial property that $u_i^+=u_i^-$ for odd $i$. 
It is unclear how it may be proven directly, 
so our proof of \eqref{u.pm} is based on a different approach.

\proof
The function
\begin{eqnarray*}
f(u) &:=& u\sqrt{g(u)/u^2},\quad u\in[-1,1],
\end{eqnarray*}
is invertible. Here the function $g(u)/u^2$ possesses the asymptotic expansion 
\begin{eqnarray*}
g(u)/u^2 &=& g_2+\sum_{i=1}^{2r} g_{i+2} u^i+o(u^{2r})
\quad\mbox{as }u\to 0.
\end{eqnarray*}
Therefore, 
\begin{eqnarray*}
f(u) &=& \sqrt{g_2}\biggl(u+\sum_{i=2}^{2r+1}f_i u^i+o(u^{2r+1})\biggr)
\quad\mbox{as }u\to 0,
\end{eqnarray*}
where $f_i$ is a polynomial of $g_3/g_2$, \ldots, $g_{i+1}/g_2$.
Let $f^{-1}(t)$ be $f$ inverse. It follows that the inverse 
function possesses an asymptotic expansion at zero up to order $2r+1$:
\begin{eqnarray}\label{expansion.f}
f^{-1}(t) &=& t/\sqrt{g_2}+\sum_{i=2}^{2r+1}c_i (t/\sqrt{g_2})^i+o(t^{2r+1})
\quad\mbox{as }t\to 0;
\end{eqnarray}
here $c_i$ is a polynomial of $g_3/g_2$, \ldots, $g_{i+1}/g_2$; 
the remainder term $o(t^{2r+1})$ may be bounded via the remainder term
in the asymptotic expansion of $g$ and the coefficients in it. 
Since the function $W(u):=\int_0^u w(u')du'$ possesses an asymptotic expansion 
at the origin up to order $2r+1$, 
\begin{eqnarray*}
W(f^{-1}(t)) &=& w_0t/\sqrt{g_2}
+\sum_{i=2}^{2r+1}\widetilde c_i (t/\sqrt{g_2})^i+o(t^r) \quad\mbox{as }t\to 0,
\end{eqnarray*}
where $\widetilde c_i$ is a polynomial of the coefficients 
$g_3/g_2$, \ldots, $g_{i+1}/g_2$, $w_0$, \ldots, $w_{i-1}$; 
here the remainder term $o(t^{2r+1})$ may be bounded via the remainder term
in the asymptotic expansions of $g$ and $W$ and the coefficients there. 
Therefore, as $t\to 0$,
\begin{eqnarray*}
W(f^{-1}(t))-W(f^{-1}(-t)) &=& 2w_0t/\sqrt{g_2}+
\sum_{i=1}^r 2\widetilde c_{2i+1} (t/\sqrt{g_2})^{2i+1}+o(t^{2r+1}).
\end{eqnarray*}
Taking into account that $u^+(t)=f^{-1}(\sqrt t)$ and
$u^-(t)=f^{-1}(-\sqrt t)$ we conclude the desired asymptotic expansion \eqref{u.pm}.
The uniform version of it---\eqref{u.pm.uni}---follows by noting that
the coefficient $g_2(\theta)$ is bounded away from zero
and that all the coefficients $\widetilde c_i(\theta)$ are polynomials 
of bounded coefficients $g_3(\theta)/g_2(\theta)$, \ldots, $g_{2r+2}(\theta)/g_2(\theta)$,
$w_0(\theta)$, \ldots, $w_{2r}(\theta)$.

Concerning \eqref{u.pm.abs}, denote
\begin{eqnarray*}
\widetilde W(u) &=& \int_0^u u' w(u')du',
\end{eqnarray*}
then
\begin{eqnarray*}
\int_{u^-(t)}^{u^+(t)} |u'|w(u')du' &=& 
\int_0^{u^+(t)} u' w(u')du'-\int_{u^-(t)}^0 u' w(u')du'\\
&=& \widetilde W(u^+(t))+\widetilde W(u^-(t)).
\end{eqnarray*}
Since the function $\widetilde W(u)$ possesses the asymptotic expansion
\begin{eqnarray*}
\widetilde W(u) &=& \frac{w_0}{2}u^2+\sum_{i=1}^{2r} \frac{w_i}{i+2} u^{i+2}+o(u^{2r+2})
\quad\mbox{as }u\to 0,
\end{eqnarray*}
we conclude from \eqref{expansion.f} that
\begin{eqnarray*}
\widetilde W(f^{-1}(t)) &=& \frac{w_0}{2g_2}t^2+
\sum_{i=3}^{2r+2} \widetilde{\widetilde c}_i (t/\sqrt{g_2})^i+o(t^{2r+2}).
\end{eqnarray*}
Therefore, as $t\to 0$,
\begin{eqnarray*}
\widetilde W(f^{-1}(t))+\widetilde W(f^{-1}(-t)) 
&=& \frac{w_0}{g_2}t^2+\sum_{i=2}^{r+1} 
2\widetilde{\widetilde c}_{2i} (t/\sqrt{g_2})^{2i}+o(t^{2r+2}).
\end{eqnarray*}
Substituting the equalities $u^+(t)=f^{-1}(\sqrt t)$ 
and $u^-(t)=f^{-1}(-\sqrt t)$, we deduce the desired 
asymptotic expansion \eqref{u.pm.abs}.

\begin{lemma}\label{l:1.f-d}
Let a function $g(\uu):\mathbb B_d\to\R^+$ possess an asymptotic expansion
$$
g(\uu)=(G_2\uu,\uu)+
\sum_{i=3}^{2r+2} g_i(\uu)+o(\|\uu\|^{2r+2}) \quad\mbox{as }\uu\to\boldsymbol 0,
$$
where $G_2$ is a positive definite matrix and 
$g_i(\uu)$ is a homogeneous polynomial of degree $i$. 
For $t>0$, denote by $B(t)$ the set of all $\uu$ such that $g(\uu)\le t$. 
Let a function $w(\uu)$ possess an asymptotic expansion
$$
w(\uu)=\sum_{i=0}^{2r} w_i(\uu)+o(\|\uu\|^{2r}) 
\quad\mbox{as }\uu\to\boldsymbol 0,
$$
where $w_i(\uu)$ is a homogeneous polynomial of degree $i$. Then
\begin{eqnarray}\label{u.Bt}
\int_{B(t)}w(\uu) d\uu &=& \frac{w_0{\rm Vol}\,\mathbb B_d}{\sqrt{\det G_2}} t^{d/2}
+\sum_{i=1}^r u_i t^{d/2+i}+o(t^{d/2+r})
\quad\mbox{as }t\downarrow 0,
\end{eqnarray}
where coefficients $u_1$, \ldots, $u_r$ only depend on 
coefficients of the polynomials $g_i$'s and $w_i$'s.
\end{lemma}

It is questionable how to extend the previous proof for multidimensional 
case $d\ge 2$. For example, if $g(\uu)=\|\uu\|^2+u_1^3$, 
then one may think of considering the invertible function
$$
f(\uu)=\uu\sqrt{g(\uu)/\|\uu\|^2}=\uu\sqrt{1+u_1^3/\|\uu\|^2}.
$$
Clearly, the function $u_1^3/\|\uu\|^2$ doesn't possesses 
an asymptotic expansion with respect to $\uu$ and 
this observation blocks the proof available in dimension $1$. 
By this reason we proceed in a different way,
by passing to hyperspherical coordinates which allows to
reduce the problem to the case $d=1$.

\proof
For $d=1$, ${\rm Vol}\,\mathbb B_1=2$ 
and the assertion is proven in Lemma \ref{l:1.f-1}.

Consider the case $d=2$. 
For $\theta\in[0,\pi)$, let $l(\theta)$ be the line passing 
through the points $(0,0)$ and $(\cos\theta,\sin\theta)$.
Since $G_2>0$, there exists a $t_0>0$ 
such that, for all $t\le t_0$ and $\theta\in[0,\pi)$, 
the set $B(t)\cap l(\theta)$ represents a segment, 
say $[b^-(\theta,t),b^+(\theta,t)]$ where $b^+(\theta,t)\in\R\times\R^+$.
Denote $u^\pm(\theta,t):=\|b^\pm(\theta,t)\|$.
Then passing to the spherical coordinates $(\theta,u)$ we deduce the equality
\begin{eqnarray}\label{int.equal}
\int_{B(t)}w(\uu) d\uu &=& 
\int_0^\pi d\theta\int_{u^-(\theta,t)}^{u^+(\theta,t)} 
|u| w(u\cos\theta,u\sin\theta)du.
\end{eqnarray}
We have
\begin{eqnarray*}
w(u\cos\theta,u\sin\theta)
&=& w_0+\sum_{i=1}^{2r} w_i(\cos\theta,\sin\theta)u^i+o(u^{2r}),
\end{eqnarray*}
so that Lemma \ref{l:1.f-1} is applicable and we conclude that
\begin{eqnarray*}
\int_{u^-(\theta,t)}^{u^+(\theta,t)} |u| w(u\cos\theta,u\sin\theta)du
&=& \frac{w_0}{g_2(\theta)} t+
\sum_{i=2}^{r+1} u_i(\theta) t^i+o(t^{r+1})
\end{eqnarray*}
as $t\downarrow 0$ uniformly on $[0,\pi)$ where 
$$
g_2(\theta)=(G_2(\cos\theta,\sin\theta),(\cos\theta,\sin\theta)).
$$
Taking into account that
$$
\int_0^\pi\frac{d\theta}{(G_2(\cos\theta,\sin\theta),(\cos\theta,\sin\theta))}
=\frac{\pi}{\sqrt{\det G_2}}
$$
and ${\rm Vol}\,\mathbb B_2=\pi$, we finally come to \eqref{u.Bt} for $d=2$.

In the same way we may proceed with an arbitrary $d\ge 3$.
First we pass to the hyperspherical coordinates
\begin{eqnarray*}
\int_{B(t)}w(\uu) d\uu &=& 
\int_{[0,\pi)^{d-1}} \det J(1,\thetaa) d\thetaa\int_{u^-(\thetaa,t)}^{u^+(\thetaa,t)} 
|u|^{d-1} w(u,\thetaa)du,
\end{eqnarray*}
then integration along the radius which is covered by Lemma \ref{l:1.f-1}.
Final integration with respect to angles completes the proof
of the asymptotic expansion and the lemma follows.

We also need the following version of Watson's lemma.

\begin{lemma}\label{multiplication}
Fix $\gamma\in\R$ and positive $y_0$, $\beta$, $c$ and $\delta$.
Then, for any $r>0$, the integral
\begin{eqnarray*}
I(x): &=& \int_{x/y_0}^\infty y^\gamma e^{-cy^\beta} (y_0-x/y)^\delta{\rm d}y
\end{eqnarray*}
possesses the expansion
\begin{eqnarray*}
I(x) &=& \Bigl(\frac{x}{y_0}\Bigr)^{1+\gamma-(1+\delta)\beta}e^{-c(x/y_0)^\beta}
\biggl(\sum_{i=0}^r I_i x^{-\beta i}+O(x^{-\beta(r+1)})\biggr)
\quad\mbox{as }x\to\infty,
\end{eqnarray*}
where $I_0=\Gamma(1+\delta)y_0^\delta/(c\beta)^{1+\delta}$.
\end{lemma}

\proof
Denote $\lambda:=c(x/y_0)^\beta$.
Changing variable $y:=(x/y_0)z^{1/\beta}$ we find that
\begin{eqnarray*}
I(x) &=& \Bigl(\frac{x}{y_0}\Bigr)^{1+\gamma}\frac{y_0^\delta}{\beta}
\int_1^\infty z^{\frac{\gamma-\delta+1}{\beta}-1}
(z^{1/\beta}-1)^\delta e^{-\lambda z}{\rm d}z\\
&=& \Bigl(\frac{x}{y_0}\Bigr)^{1+\gamma}\frac{y_0^\delta}{\beta}
e^{-\lambda}
\int_0^\infty (1+u)^{\frac{\gamma+1}{\beta}-1}
(1-(1+u)^{-1/\beta})^\delta e^{-\lambda u}{\rm d}u.
\end{eqnarray*}
If $r>\frac{\gamma+1}{\beta}-2-\delta$ then, for all $u>0$,
$$
\biggl|(1+u)^{\frac{\gamma+1}{\beta}-1}
\biggl(\frac{1-(1+u)^{-1/\beta}}{u}\biggr)^\delta
-\frac{1}{\beta^\delta}+\sum_{i=1}^r c'_i u^i\biggr| \le c'u^{r+1}
$$
for some $c_i'$ and $c'<\infty$. Hence,
$$
\biggl|I(x)-\Bigl(\frac{x}{y_0}\Bigr)^{1+\gamma} e^{-\lambda}
\sum_{i=0}^r I'_i \int_0^\infty u^{\delta+i}e^{-\lambda u}{\rm d}u\biggr| 
\le I'x^{1+\gamma}e^{-\lambda}\int_0^\infty u^{\delta+r+1}e^{-\lambda u}{\rm d}u,
$$
where $I'_0=y_0^\delta/\beta^{1+\delta}$ and 
$I'_i=c_i'y_0^\delta/\beta$ for $i\ge 1$. In its turn,
\begin{eqnarray*}
\int_0^\infty u^{\delta+i}e^{-\lambda u}{\rm d}u
&=& \frac{\Gamma(\delta+i+1)}{\lambda^{1+\delta+i}}
=\frac{\Gamma(\delta+i+1)}{c^{1+\delta+i}}
\Bigl(\frac{x}{y_0}\Bigr)^{-\beta(1+\delta+i)},
\end{eqnarray*}
which completes the proof.

Let us proceed with the proof of Theorem \ref{th:m=0}.
Substituting the result of Lemma \ref{asy.for.g.mu.0}
into \eqref{multipl.our} we deduce that, as $x\to\infty$
\begin{eqnarray*}
\Pr\{g(\etaa)>x\} &\sim& \frac{1}{\sqrt{2\pi}\alpha\Gamma((d+1)/2)}
\sum_{j=1}^k \frac{|\det J(1,\phii_j)|}{\sqrt{|\det g''(\phii_j)|}}\\
&&\hspace{40mm} \times \int_{x/\hat g}^\infty y^{d/\alpha-1}e^{-y^{2/\alpha}/2}
(\hat g-x/y)^{(d-1)/2}{\rm d}y.
\end{eqnarray*}
Now we apply Lemma \ref{multiplication} with
$y_0:=\hat g$, $\gamma:=d/\alpha-1$, $\beta:=2/\alpha$, $c:=1/2$
and $\delta:=(d-1)/2$ we deduce from that the following
asymptotics, as $x\to\infty$:
\begin{eqnarray*}
\Pr\{g(\etaa)>x\} &\sim&
\frac{(\alpha\hat g)^{\frac{d-1}{2}}}{\sqrt{2\pi}}
\sum_{j=1}^k \frac{|\det J(1,\phii_j)|}{\sqrt{|\det g''(\phii_j)|}}
\Bigl(\frac{x}{\hat g}\Bigr)^{-1/\alpha} e^{-(x/\hat g)^{2/\alpha}/2},
\end{eqnarray*}
which completes the proof of the tail asymptotics.

Next, we  prove the claim in \eqref{eqden} which shows a tractable 
expression of the density of $g(\etaa)$ in terms of tail 
characteristics of $g(\etaa)$.

\begin{lemma}\label{dens.via.tail}
The density function $p_{g(\eta)}$ of distribution of
$g(\etaa)$ restricted to $\R\setminus\{0\}$ exists
and possesses the representation \eqref{eqden}. Moreover,
\begin{eqnarray}\label{eqden.2}
p_{g(\etaa)}(x) &=&
\frac{1}{\alpha x}\biggl( \frac{1}{\alpha2^{d/2-1}\Gamma(d/2)}
\int_{x/\hat g}^\infty y^{(d+2)/\alpha-1}e^{-y^{2/\alpha}/2}
\Pr\{g(\zetaa)>x/y\}{\rm d}y\nonumber\\
&&\hspace{65mm}-d\cdot\Pr\{g(\etaa)>x\}\biggr),
\end{eqnarray}
where $\zetaa$ is uniformly distributed on $\Sp_{d-1}$.
\end{lemma}

\proof
Since $\etaa$ is a standard Gaussian random vector,
\begin{eqnarray*}
p_{g(\etaa)}(x) &=& -\frac{{\rm d}}{{\rm d}x}\Pr\{g(x^{-1/\alpha}\etaa)>1\}\\
&=& -\frac{{\rm d}}{{\rm d}x}\biggl(\frac{x^{d/\alpha}}{(2\pi)^{d/2}}
\int_{\{\boldsymbol v\in\R^d:g(\boldsymbol v)>1\}}
e^{-x^{2/\alpha}\|\boldsymbol v\|^2/2}{\rm d}\boldsymbol v\biggr)\\
&=& -\frac{d}{\alpha x}\frac{x^{d/\alpha}}{(2\pi)^{d/2}}
\int_{\{\boldsymbol v\in\R^d:g(\boldsymbol v)>1\}}
e^{-x^{2/\alpha}\|\boldsymbol v\|^2/2}{\rm d}\boldsymbol v\\
&&\hspace{30mm}+\frac{x^{2/\alpha-1}}{\alpha}
\frac{x^{d/\alpha}}{(2\pi)^{d/2}}
\int_{\{\boldsymbol v\in\R^d:g(\boldsymbol v)>1\}}
|\boldsymbol v|^2
e^{-x^{2/\alpha}\|\boldsymbol v\|^2/2}{\rm d}\boldsymbol v.
\end{eqnarray*}
Therefore,
\begin{eqnarray*}
p_{g(\etaa)}(x) &=& -\frac{d}{\alpha x}\Pr\{g(x^{-1/\alpha}\etaa)>1\}
+\frac{x^{2/\alpha-1}}{\alpha}\E\{\|x^{-1/\alpha}\etaa\|^2;
g(x^{-1/\alpha}\etaa)>1\},
\end{eqnarray*}
which implies \eqref{eqden}.

Similarly to \eqref{multipl.our}, 
we derive from \eqref{chi.alpha} the equality
\begin{eqnarray*}
\E\{\|\etaa\|^2;g(\etaa)>x\} &=&
\E\bigl\{\E\bigl\{(\chi^\alpha)^{2/\alpha}{\mathbb I}\{g(\zetaa)>x/\chi^\alpha\}
\bigr\}\mid\chi\bigr\}\\
&=& \int_{x/\hat g}^\infty y^{2/\alpha}
p_{\chi^\alpha}(y)\Pr\{g(\zetaa)>x/y\}{\rm d}y\\
&=& \frac{1}{\alpha2^{d/2-1}\Gamma(d/2)}
\int_{x/\hat g}^\infty y^{(d+2)/\alpha-1}e^{-y^{2/\alpha}/2}
\Pr\{g(\zetaa)>x/y\}{\rm d}y
\end{eqnarray*}
establishing thus the proof. 

Further application of Lemmas \ref{asy.for.g.mu.0} and
\ref{multiplication} completes the proof of the density
function asymptotic expansion.

In the Gaussian case, yet another approach for estimating the tail
of Gaussian chaos seems to be applicable. Consider $n$ independent
copies $\etaa_1$, \ldots, $\etaa_n$ of $\etaa$, then
$$
\Pr\{g(\etaa)>x\} = \Pr\Bigl\{g\Bigl(\frac{\etaa_1 +\ldots+\etaa_n}{\sqrt n}\Bigr)>x\Bigr\}
= \Pr\Bigl\{g\Bigl(\frac{\etaa_1 +\ldots+\etaa_n}{n}\Bigr)
>\frac{x}{n^{\alpha/2}}\Bigr\}.
$$
Therefore, considering $x=tn^{\alpha/2}$, we have
$$
\Pr\{g(\etaa)>x\} = \Pr\Bigl\{g\Bigl(\frac{\etaa_1 +\ldots+\etaa_n}{n}\Bigr)>t\Bigr\}.
$$
Hence, this reduces the problem of the tail behavior as $x\to\infty$
to that of large deviation as $n\to\infty$. Then one may try 
to apply some results on asymptotic expansions in large deviations,
for the distribution as well as for the density, see e.g., Borovkov
and Rogozin \cite{BR1965}. We just mention that it is easily seen 
that the integration over a domain of the asymptotic expansion for 
the large deviation probabilities is not simpler than our integration 
related to a chi-squared distribution. 

We conclude this section with the proof of the equality
\eqref{cart.local} which follows from the following result.

\begin{lemma}\label{l:car.local}
Under the  conditions of Theorem \ref{th:m=0} at every point
$\boldsymbol v_j=h_j(\boldsymbol z_j)\in\mathcal M$,
$\boldsymbol z_j\in (0,2)^{d-1}\times\{1\}$
\begin{eqnarray*}
(J_j^{-1}(\boldsymbol z_j))^{\rm T}(g\circ h_j)_{d-1}''(\boldsymbol z_j)
J_j^{-1}(\boldsymbol z_j)
&=& g''_{d-1}(\boldsymbol v_j)-\hat g\alpha I_{d-1}.
\end{eqnarray*}
\end{lemma}

\proof
Indeed, since $\boldsymbol v_j$ is the point of the maximum of
the function $g$, Taylor's expansion at this points reads as follows:
with $\boldsymbol v=h_j(\boldsymbol z)$,
\begin{eqnarray*}
g(\boldsymbol v) &=& \hat g+\frac{1}{2}\bigl((g\circ h_j)_{d-1}''(\boldsymbol z_j)
(\boldsymbol z-\boldsymbol z_j),\boldsymbol z-\boldsymbol z_j\bigr)
+o(\|\boldsymbol v-\boldsymbol v_j\|^2)\\
&=&  \hat g+\frac{1}{2}\bigl((g\circ h_j)_{d-1}''(\boldsymbol z_j)
(h_j^{-1}(\boldsymbol v)-h_j^{-1}(\boldsymbol v_j)),
h_j^{-1}(\boldsymbol v)-h_j^{-1}(\boldsymbol v_j)\bigr)
+o(\|\boldsymbol v-\boldsymbol v_j\|^2)
\end{eqnarray*}
as $\boldsymbol v\to\boldsymbol v_j$. Since
$$
h_j^{-1}(\boldsymbol v)-h_j^{-1}(\boldsymbol v_j)
=J_j^{-1}(\boldsymbol z_j)(\boldsymbol v-\boldsymbol v_j)
+o(\|\boldsymbol v-\boldsymbol v_j\|),
$$
we obtain that
\begin{eqnarray}\label{hes.0}
g(\boldsymbol v) &=&
\hat g+\frac{1}{2}\bigl((g\circ h_j)''(\boldsymbol z_j)
J_j^{-1}(\boldsymbol z_j)(\boldsymbol v-\boldsymbol v_j),
J_j^{-1}(\boldsymbol z_j)(\boldsymbol v-\boldsymbol v_j)\bigr)
+o(\|\boldsymbol v-\boldsymbol v_j\|^2)\nonumber\\[-2mm]
\end{eqnarray}
as $\boldsymbol v\to\boldsymbol v_j$.
Consider the projection $\boldsymbol u$ of the point
$\boldsymbol v$ onto the hyperplane
$(u_1-(\boldsymbol v_j)_1)(\boldsymbol v_j)_1+\ldots
+(u_d-(\boldsymbol v_j)_d)(\boldsymbol v_j)_d=0$.
The equalities
\begin{eqnarray*}
\|\boldsymbol v-\boldsymbol u\| &=&
1-(\boldsymbol v,\boldsymbol v_j)
=(\boldsymbol v,\boldsymbol v-\boldsymbol v_j)\\
&=& (\boldsymbol v-\boldsymbol v_j+\boldsymbol v_j,\boldsymbol v-\boldsymbol v_j)\\
&=& \|\boldsymbol v-\boldsymbol v_j\|^2
+(\boldsymbol v_j,\boldsymbol v)-1
\end{eqnarray*}
yield that
\begin{eqnarray}\label{v-u}
\|\boldsymbol v-\boldsymbol u\| &=&
\|\boldsymbol v-\boldsymbol v_j\|^2/2
=\|\boldsymbol u-\boldsymbol v_j\|^2/2
+o(\|\boldsymbol u-\boldsymbol v_j\|^4).
\end{eqnarray}
Applying this in \eqref{hes.0} we deduce that
\begin{eqnarray}\label{hes.1}
g(\boldsymbol v) &=& \hat g+\frac{1}{2}
\bigl((J_j^{-1}(\boldsymbol z_j))^{\rm T}(g\circ h_j)_{d-1}''(\boldsymbol z_j)
J_j^{-1}(\boldsymbol z_j)
(\boldsymbol u-\boldsymbol v_j), \boldsymbol u-\boldsymbol v_j\bigr)
+o(\|\boldsymbol u-\boldsymbol v_j\|^2).\nonumber\\[-2mm]
\end{eqnarray}
On the other hand, again by Taylor's expansion
\begin{eqnarray}\label{hes.2}
g(\boldsymbol u) &=& \hat g+\frac{1}{2}(g''_{d-1}(\boldsymbol v_j)
(\boldsymbol u-\boldsymbol v_j),\boldsymbol u-\boldsymbol v_j)
+o(\|\boldsymbol u-\boldsymbol v_j\|^2)
\quad\mbox{ as }\boldsymbol u\to\boldsymbol v_j.
\end{eqnarray}
In addition
\begin{eqnarray}\label{hes.3}
g(\boldsymbol v) &=& g(\boldsymbol u)
+(\nabla g(\boldsymbol u),\boldsymbol v-\boldsymbol u)
+O(\|\boldsymbol v-\boldsymbol u\|^2)\nonumber\\
&=& g(\boldsymbol u)
+(\nabla g(\boldsymbol v_j),\boldsymbol v-\boldsymbol u)
+O(\|\boldsymbol u-\boldsymbol v_j\|^2)
\end{eqnarray}
because $\nabla g(\boldsymbol u)\to\nabla g(\boldsymbol v_j)$
and $\boldsymbol v-\boldsymbol u=O(\|\boldsymbol u-\boldsymbol v_j\|^2)$.
Since $\boldsymbol v_j$ and $\boldsymbol v-\boldsymbol u$ are collinear,
\begin{eqnarray*}
(\nabla g(\boldsymbol v_j),\boldsymbol v-\boldsymbol u) &=&
-\|\boldsymbol v-\boldsymbol u\|\lim_{\varepsilon\to 0}
\frac{g(\boldsymbol v_j+\varepsilon\boldsymbol v_j)-g(\boldsymbol v_j)}{\varepsilon}.
\end{eqnarray*}
By the homogeneity of the function $g$,
\begin{eqnarray*}
\lim_{\varepsilon\to 0}
\frac{g(\boldsymbol v_j+\varepsilon\boldsymbol v_j)-g(\boldsymbol v_j)}{\varepsilon}
&=& g(\boldsymbol v_j)\lim_{\varepsilon\to 0}
\frac{(1+\varepsilon)^\alpha-1}{\varepsilon} = -\hat g\alpha,
\end{eqnarray*}
so that we have
\begin{eqnarray*}
(\nabla g(\boldsymbol v_j),\boldsymbol v-\boldsymbol u)
&=& -\hat g\alpha\|\boldsymbol v-\boldsymbol u\|\\
&=& -\frac{1}{2}\hat g\alpha
(\boldsymbol u-\boldsymbol v_j,\boldsymbol u-\boldsymbol v_j)
+o(\|\boldsymbol u-\boldsymbol v_j\|^4),
\end{eqnarray*}
by the equality \eqref{v-u}.
Combining \eqref{hes.1}, \eqref{hes.2} and \eqref{hes.3} we conclude the desired equality of the matrices.

\section{Proof of Theorem \ref{th:m=1.d-2}}
\label{sec:proof.d}

Since $\mathcal M$ is $C^{2r+2}$-smooth manifold in $\Sp_{d-1}$, 
there exists some neighborhood $U$ of $\mathcal M$ in $\R^d$
such that it may be partitioned into a finite number of disjoint
sets $U_1$, \ldots, $U_n$ such that, for every $1\le j\le n$,
the manifold $\mathcal M\cap U_j$ is elementary,
that is, there exists some bijection $h_j:[0,2]^d\to{\rm cl}\,(U_j)$
(the closure of $U_j$)
which is $2r+2$ times differentiable, non-degenerate and such that
$$
h_j([0,2]^{d-1}\times\{1\})=\Sp_{d-1}\cap{\rm cl}\,(U_j)
\quad\mbox{and}\quad
h_j([0,2]^m\times\{1\}^{d-m})=\mathcal M\cap{\rm cl}\,(U_j).
$$
It is non-degenerate in the sense that its Hessian
is non-zero at every point $\boldsymbol z\in[0,2]^d$.

The proof of Theorem \ref{th:m=1.d-2} follows the lines of the previous one.
The main difference consists in the estimation of the
probability $\Pr\{g(\zetaa)>\hat g-t\}$ for small values of $t>0$.
Because of this, we only need to show the following result.

\begin{lemma}\label{asy.for.g.mu}
Under the  conditions of Theorem \ref{th:m=1.d-2}
$$
\Pr\{g(\zetaa)>\hat g-t\} =
\sum_{i=0}^r g_i t^{\frac{d-1-m}{2}+i}
+o(t^{\frac{d-1-m}{2}+r})\quad\mbox{as }t\downarrow 0,
$$
where
$$
g_0=2^{\frac{d-1-m}{2}}\frac{{\rm Vol}\,\mathbb B_{d-1-m}}{{\rm mes}\,\Sp_{d-1}}
\int_{\mathcal M}
\bigl|\det\bigl(g_{d-1-m}''(\boldsymbol v)
-(\alpha\hat g)I_{d-1-m}\bigr)\bigr|^{-1/2}{\rm d}V.
$$
\end{lemma}

Since ${\rm Vol}\,\mathbb B_{d-1-m}=\frac{\pi^{(d-1-m)/2}}{\Gamma((d+1-m)/2)}$
and ${\rm mes}\,\mathbb S_{d-1}=\frac{2\pi^{d/2}}{\Gamma(d/2)}$,
we have the following alternative representation for the
constant $g_0$, in terms of the hyperspherical coordinates:
$$
g_0=\frac{2^{d/2-1}}{(2\pi)^{(1+m)/2}}
\frac{\Gamma(d/2)}{\Gamma((d+1-m)/2)} \int_{\mathcal M_\varphi}
\frac{|\det J(1,\phii)|}{\sqrt{|\det g_{d-1-m}''(\phii)|}}{\rm d}\Vfi.
$$

\proof
For every $j\le n$, consider a random vector $\nuu_j$ valued in
$[0,2]^{d-1}$ with density function
$\frac{|\det J_j(\boldsymbol z)|}{{\rm mes}\,(\Sp_{d-1}\cap U_j)}$,
$\boldsymbol z\in[0,2]^{d-1}$, where $J_j(\boldsymbol z)$
is the Jacobian matrix of the function $h_j$ restricted to the first
$d-1$ coordinates. Then $h_j(\nuu_j,1)$ has the uniform distribution
on the set $\Sp_{d-1}\cap U_j$.

Let $t>0$ be so small that the $t$-neighborhood of $\mathcal M$
is contained in the set $U$.
Consider the following decomposition:
\begin{eqnarray*}
\Pr\{g(\zetaa)>\hat g-t\} &=&
\sum_{j=1}^n\Pr\{g(\zetaa)>\hat g-t,\zetaa\in\Sp_{d-1}\cap U_j\}\\
&=& \sum_{j=1}^n\Pr\{g(h_j(\nuu_j,1))>\hat g-t\}
\frac{{\rm mes}\,(\Sp_{d-1}\cap U_j)}{{\rm mes}\,\Sp_{d-1}}
\end{eqnarray*}
and compute the asymptotic behaviour of the $j$th term on the right.
Since $h_j([0,2]^m\times\{1\}^{d-m})=\mathcal M\cap{\rm cl}\,(U_j)$,
we have $g(h_j(\boldsymbol s,1,\ldots,1))=\hat g$ for every
point $\boldsymbol s\in[0,2]^m$. The function
$g(h_j(\boldsymbol s,\cdot,\ldots,\cdot,1))$ of $d-1-m$ arguments
is $2r+2$ times differentiable. Then the same arguments as in
Lemma \ref{asy.for.g.mu.0} yield the decomposition
\begin{eqnarray*}
\Pr\{g(h_j(\nuu_j,1))>\hat g-t\mid
\nuu_j\in\{\boldsymbol s\}\times[0,2]^{d-1-m}\} &=&
\sum_{i=0}^r g_{ji}(\boldsymbol s) t^{\frac{d-1-m}{2}+i}
+o(t^{\frac{d-1-m}{2}+r})
\end{eqnarray*}
as $t\downarrow 0$ where
$$
g_{j0}(\boldsymbol s)=
2^{\frac{d-1-m}{2}}\frac{{\rm Vol}\,\mathbb B_{d-1-m}}
{{\rm mes}\,(\Sp_{d-1}\cap U_j)}
\frac{|\det J_j(\boldsymbol s,1,\ldots,1)|}
{\sqrt{|\det (g\circ h_j)''(\boldsymbol s,1,\ldots,1)|}},
$$
where the Hessian of $g\circ h_j$ is taken with respect to the
last $d-1-m$ arguments. Integration over $\boldsymbol s\in[0,2]^m$
finally implies that
\begin{eqnarray*}
\Pr\{g(h_j(\nuu_j,1))>\hat g-t\} &=&
\sum_{i=0}^r g_{ji} t^{\frac{d-1-m}{2}+i}
+o(t^{\frac{d-1-m}{2}+r})
\end{eqnarray*}
as $t\downarrow 0$ where
$$
g_{j0}=2^{\frac{d-1-m}{2}}\frac{{\rm Vol}\,\mathbb B_{d-1-m}}
{{\rm mes}\,(\Sp_{d-1}\cap U_j)} \int_{[0,2]^m}
\frac{|\det J_j(\boldsymbol s,1,\ldots,1)|}
{\sqrt{|\det (g\circ h_j)''(\boldsymbol s,1,\ldots,1)|}}
{\rm d}\boldsymbol s,
$$
which proves the lemma.

\section{Proof of Theorem \ref{thPolar}}
\label{sec:proof.3} 

The crucial step of the proof is again to find the tail 
asymptotics of $g(\zetaa)$.
As in the proof of Lemma \ref{asy.for.g.mu} we have
$$
\Pr\{g(\zetaa)>\hat g-t\} \sim
g_0 t^{\frac{d-1-m}{2}} \quad\mbox{as }t\downarrow 0,
$$
where
\begin{eqnarray}
g_0=\frac{(2\pi)^{\frac{d-1-m}{2}}}{\Gamma((d+1-m)/2)} \int_{\mathcal M_\varphi}
\frac{p_\nuu(\phii)}{\sqrt{|\det g_{d-1-m}''(\phii)|}}{\rm d}\Vfi.
\end{eqnarray}
Next---here we follow a simplified version compared to Hashorva (2012)---
\begin{eqnarray}\label{pre.asy.g}
\Pr\{g(\etaa)>x\} &=& \Pr\{\chi^\alpha g(\zetaa)>x\}\nonumber\\
&=& \int_0^\infty \Pr\Bigl\{g(\zetaa)>\frac{x}{x/\hat g+y}\Bigr\}
\Pr\{\chi^\alpha\in x/\hat g+{\rm d}y\}\nonumber\\
&=& \int_0^\infty \Pr\Bigl\{g(\zetaa)>\frac{x}{x/\hat g+y/w(x/\hat g)}\Bigr\}
\Pr\Bigl\{\chi^\alpha\in \frac{x}{\hat g}+\frac{{\rm d}y}{w(x/\hat g)}\Bigr\}\nonumber\\
&=& \Pr\Bigl\{\chi^\alpha > \frac{x}{\hat g}\Bigr\}
\int_0^\infty \Pr\Bigl\{g(\zetaa)>\hat g-\frac{\hat g^2}{xw(x/\hat g)}
\frac{y}{1+y\hat g/xw(x/\hat g)}\Bigr\}\nonumber\\
&&\hspace{40mm} \Pr\Bigl\{\chi^\alpha\in \frac{x}{\hat g}+\frac{{\rm d}y}{w(x/\hat g)}
\Big|\chi^\alpha > \frac{x}{\hat g}\Bigr\}.
\end{eqnarray}
Since $xw(x)\to\infty$ as $x\to\infty$, 
\begin{eqnarray*}
\Pr\Bigl\{g(\zetaa)>\hat g-\frac{\hat g^2}{xw(x/\hat g)}
\frac{y}{1+y\hat g/xw(x/\hat g)}\Bigr\}
&\sim& g_0 \Bigl(\frac{\hat g^2y}{xw(x/\hat g)}\Bigr)^{\frac{d-1-m}{2}}
\end{eqnarray*}
as $x\to\infty$ uniformly on any $y$-compact set. In addition,
\begin{eqnarray*}
\Pr\Bigl\{g(\zetaa)>\hat g-\frac{\hat g^2}{xw(x/\hat g)}
\frac{y}{1+y\hat g/xw(x/\hat g)}\Bigr\}
&\le& \Pr\Bigl\{g(\zetaa)>\hat g-\frac{y\hat g^2}{xw(x/\hat g)}\Bigr\}\\
&\le& c \Bigl(\frac{y}{xw(x/\hat g)}\Bigr)^{\frac{d-1-m}{2}},
\end{eqnarray*}
for some $c<\infty$. The latter asymptotics and upper bound allow
to apply Lebesgue's dominated convergence theorem in \eqref{pre.asy.g}
and to conclude that, as $x\to\infty$,
\begin{eqnarray*}
\lefteqn{\Pr\{g(\etaa)>x\}}\\
 &\sim& g_0 \Bigl(\frac{\hat g^2}{xw(x/\hat g)}\Bigr)^{\frac{d-1-m}{2}}
\Pr\Bigl\{\chi^\alpha > \frac{x}{\hat g}\Bigr\}
\int_0^\infty y^{\frac{d-1-m}{2}}
\Pr\Bigl\{\chi^\alpha\in \frac{x}{\hat g}+\frac{{\rm d}y}{w(x/\hat g)}
\Big|\chi^\alpha > \frac{x}{\hat g}\Bigr\}.
\end{eqnarray*}
It follows from the Davis--Resnick tail property---see 
\cite[Proposition 1.1]{DR}---that, for any fixed $\gamma>0$, 
there exists a $c_1<\infty$ such that for all $u$, $v>0$
\begin{eqnarray*}
\Pr\{\chi^\alpha > u+v/w(u)\mid\chi^\alpha > u\} 
&\le& c_1/v^\gamma.
\end{eqnarray*}
This ensures the following convergence of moments
\begin{eqnarray*}
\int_0^\infty y^{\frac{d-1-m}{2}}
\Pr\Bigl\{\chi^\alpha\in \frac{x}{\hat g}+\frac{{\rm d}y}{w(x/\hat g)}
\Big|\chi^\alpha > \frac{x}{\hat g}\Bigr\}
&\to& \int_0^\infty y^{\frac{d-1-m}{2}}e^{-y}{\rm d}y\\
&=& \Gamma\Bigl(\frac{d+1-m}{2}\Bigr),
\end{eqnarray*}
and hence the first claim follows.

Since further the scaling function $w(\cdot)$ 
is self-neglecting (see e.g., Resnick (1987)) i.e.,
$$
w(t+s/w(t)) \sim w(t)\quad\mbox{as }t\uparrow x_+
$$
locally uniformly in $s$, then $g(\etaa)$ is also in the Gumbel
max-domain of attraction with the same scaling function $w$ as $\chi^\alpha$. 
Thus the second claim follows.

\section*{Acknowledgment}

The authors gratefully acknowledge helpful consultation on random matrices 
by Vadim Gorin and valuable discussions with Philippe Barbe.
They are also very thankful to two referees whose comments 
helped a lot to improve the paper.

The authors kindly acknowledge partial support from SNSF grants 200021-140633/1, 200021-134785 and 
RARE--318984 (an FP7 Marie Curie IRSES Fellowship).
The research of V. I. Piterbarg is supported by the Russian Foundation 
for Basic Research, Projects 11-01-00050-a and 14-01-00075.


\begin{thebibliography}{99}

\bibitem{Are} 
Arendarczyk, M., D\c{e}bicki, K. (2011) 
Asymptotics of supremum distribution of a Gaussian process over a Weibullian time.  
\emph{Bernoulli}, {\bf 17}, 194--210.

\bibitem{AGZ}
Anderson, G. V., Guionnet, A.,  and Zeitouni, O. (2010)
{\it An Introduction to Random Matrices}.
Cambridge: Cambridge University Press.

\bibitem{Gine}
Arcones, M. A., and Gin\'e, E. (1993)
On decoupling, series expansions, and tail behavior of chaos processes.
\textit{J. Theor. Probab.} {\bf 6} 101--122.

\bibitem{BKR}
Balkema, A. A., Kl\"uppelberg, C., and Resnick, S. I. (1993)
Densities with Gaussian tails.
\textit{Proc. London Math. Soc.} {\bf 66} 568--588.

\bibitem{Barbe}
Barbe, Ph. (2003)
{\it Approximation of integrals over asymptotic sets
with applications to probability and statistics}.
arXiv:math/0312132.

\bibitem{Breiman}
Breiman, L. (1965)
On some limit theorems similar to the arc-sin law.
\textit{Theory Probab. Appl.} {\bf 10}, 323-331.

\bibitem{borel}
Borell, C. (1978)
Tail probabilities on Gauss space.
{\it Lect. Notes in Math.} {\bf 644}, Springer, Berlin, 73--82.

\bibitem{BR1965}
Borovkov, A. A., and Rogozin, B. A. (1965)
On the multi-dimensional central limit theorem. 
{\it Theory Probab. Appl.} {\bf 10} 55--62.

\bibitem{BR1996}
Breitung, K., and Richter, W.-D. (1996)
A geometric approach to an asymptotic expansion for 
large deviation probabilities of Gaussian random vectors.
{\it J. Multiv. Anal.} {\bf 58}, 1--20.

\bibitem{Cline1994}
Cline, D. B. H., and Samorodnitsky, G. (1994)
Subexponentiality of the product of independent random variables.
\textit{Stochastic Processes Appl.} {\bf 49}, 75--98.

\bibitem{DR}
Davis, R. A., and Resnick, S. I. (1988)
Extremes of moving averages of random variables from the domain
of attraction of the double exponential distribution.
{\it Stoch. Process. Appl.} {\bf 30}, 41--68.


\bibitem{EMB} 
Embrechts, P., Kl\"{u}ppelberg, C., and Mikosch, T. (1997)
\textit{Modelling Extreme Events for Insurance and Finance}. 
Springer-Verlag, Berlin.

\bibitem{Fedoryuk}
Fedoryuk, M. V. (1977)
\textit{Metod perevala (The saddlepoint method)}. 
Nauka, Moscow [In Russian].

\bibitem{FKZ}
Foss, S., Korshunov, D., and Zachary, S. (2011)
\textit{An Introduction to Heavy-Tailed and Subexponential Distributions}. Springer, New York.

\bibitem{HW}
Hanson, D. L., and Wright, F. T. (1971)
A bound on tail probabilities for quadratic forms in independent random variables.
\textit{Ann. Math. Stat.} {\bf 42}, 1079--1083.

\bibitem{HE4}
Hashorva, E. (2015) Extremes of aggregated Dirichlet risks.
{\it J. Multivariate Anal.} {\bf 133}, 334--345.

\bibitem{HE}
Hashorva, E. (2012)
Exact tail asymptotics in bivariate scale mixture models. 
{\it Extremes} {\bf 15}, 109--128.

\bibitem{HE0}
Hashorva, E. (2010) Asymptotics of the norm of elliptical random vectors.  {\it J. Multivariate Anal.} {\bf 101}, 926--935.

\bibitem{HJT}
Hashorva, E., Ji., L., and Tan, Z. (2012)
On the infinite sums of deflated Gaussian products.
{\it Elect. Comm. Probab.} {\bf 17}, 1--8.

\bibitem{HKP}
Hashorva, E., Korshunov, D. A., and Piterbarg, V. I. (2013)
On extremal behavior of Gaussian chaos. 
{\it Doklady Mathematics} {\bf 88}, 566--568.


\bibitem{hoeffding}
Hoeffding, W. (1964)
On a theorem of V.N. Zolotarev.
\textit{Theory Probab. Appl.} {\bf 9}, 89--91.

\bibitem{husler}
H\"usler, J., Liu, R., and Singh, K. (2002)
A formula for the tail probability of a multivariate
normal distribution and its applications.
{\it J. Multiv. Anal.} {\bf 82}, 422-"1¤7430.

\bibitem{imkeller}
Imkeller, P. (1994)
On exact tails for limiting distributions of $U$-statistics
in the second Gaussian chaos. \textit{Chaos expansions, multiple Wiener--It\^o integrals and
their applications (Guanajuato, 1992)}, 239--244, 
{\it Probab. Stochastics Ser., CRC, Boca Raton, FL}.

\bibitem{ivanoff}
Ivanoff, B. G., and Weber, N. C. (1998)
Tail probabilities for weighted sums of products of normal random variables.
\textit{Bull. Austral. Math. Soc.} {\bf 58}, 239--244.

\bibitem{jaconsen_etall}
Jacobsen, M., Mikosch, T., Rosinski, J., and Samorodnitsky, G. (2009)
Inverse problems for regular variation of linear filters, a cancellation
property for sigma-finite measures, and identification of stable laws.
\textit{Ann. Appl. Probab.} {\bf 19}, 210--242.

\bibitem{janson}
Janson, S. (1997) 
{\it Gaussian Hilbert Spaces}.
Cambridge University Press.

\bibitem{Major2005}
Major, P. (2005)
Tail behaviour of multiple random integrals and $U$-statistics.
{\it Probab. Surv.} {\bf 2}, 448--505.

\bibitem{Major2007}
Major, P. (2007)
On a multivariate version of Bernstein's inequality.
{\it Electron. J. Probab.} {\bf 12}, 966--988 (electronic).

\bibitem{MR}
Matthews, P. C. and Rukhin, A. L. (1993)
Asymptotic distribution of the normal sample range.
{\it Ann. Appl. Probab.}, {\bf 3}, 454--466

\bibitem{Latala1999}
Lata{\l }a, R. (1999)
Tail and moment esimates for some types of chaos.
\textit{Studia Mathematica} {\bf 135}, 39--53.

\bibitem{Latala2006}
Lata{\l }a, R. (2006)
Estimates of moments and tails of Gaussian chaoses.
\textit{Ann. Probab.} {\bf 34}, 2315--2331.

\bibitem{lehec}
Lehec, J. (2011)
Moments of the Gaussian chaos.
\textit{S\'eminaire de Probabilit\'es XLIII}, 
Lecture Notes in Math. 2006, Springer, 327--340.

\bibitem{PR1988}
Pap, G., and Richter, W.-D. (1988)
Zum asymptotischen Verhalten der Verteilungen and der Dichten
gewisser Funktionale Gauss'scher Zufallsvektoren.
{\it Math. Nachr.} {\bf 135}, 119--124.

\bibitem{Piterbarg}
Piterbarg, V. I. (1994)
High excursions for nonstationary generalized chi-square processes.
\textit{Stochastic Processes Appl.} {\bf 53}, 307--337.

\bibitem{Piterbarg1996}
Piterbarg, V. I. (1996)
{\it Asymptotic Methods in the Theory of Gaussian Processes and Fields}.
In: Transl. Math. Monographs, vol. 148.
AMS, Providence, RI.

\bibitem{Prekopa}
Pr\'{e}kopa, A. (1967) On random determinants I.  
{\it Stud. Sci. Math. Hung.} {\bf 2},  125--132.

\bibitem{Resnick}
Resnick, S. I. (1987)
{\it Extreme Values, Regular Variation and Point Processes.}
Springer, New York.

\bibitem{Rootzen}
Rootz\'en, H. (1987)
A ratio limit theorem for the tails of weighted sums.
{\it Ann. Probab.} {\bf 15}, 728--747.

\bibitem{Sornette}
Sornette, D. (1998)
Multiplicative processes and power laws.
\textit{Phys. Rev. E} {\bf 57}, 4811--4813.

\bibitem{Wiener}
Wiener, N. (1938)
The homogeneous chaos.
{\it Amer. J. Math.} {\bf 60}, 897--936.

\bibitem{zolotarev}
Zolotarev, V. M. (1961)
Concerning a certain probability problem.
\textit{Theory Probab. Appl.} {\bf 6(2)}, 201--204.

\end{thebibliography}
\end{document}